\documentclass[11pt,a4paper,leqno]{article}
\usepackage{a4wide}
\usepackage{latexsym}
\usepackage{amsmath}
\usepackage{amssymb}
\usepackage{stackrel}  
\usepackage{color} 
\usepackage{graphicx}
\usepackage[linesnumbered,ruled,vlined]{algorithm2e}

\usepackage{authblk}

\usepackage{hyperref}

\newtheorem{defin}{Definition}
\newtheorem{lemma}{Lemma}
\newtheorem{prop}{Proposition}
\newtheorem{theo}{Theorem}
\newtheorem{corol}{Corollary}

\newenvironment{proof}{\medskip\par\noindent{\bf Proof}}{\hfill $\Box$
\medskip\par}

\newcommand{\C}{\mathbb{C}}
\newcommand{\N}{\mathbb{N}}
\newcommand{\R}{\mathbb{R}}
\newcommand{\Z}{\mathbb{Z}}

\begin{document}
\title{On some $q$-analog of an initial value problem with infinite order irregular singularity and Mahler transforms}

\author[1]{Alberto Lastra}
\author[2]{St\'ephane Malek}
\affil[1]{Universidad de Alcal\'a, Dpto. F\'isica y Matem\'aticas, Alcal\'a de Henares, Madrid, Spain. {\tt alberto.lastra@uah.es}}
\affil[2]{University of Lille, Laboratoire Paul Painlev\'e, Villeneuve d'Ascq cedex, France. {\tt stephane.malek@univ-lille.fr}}


\date{}

\maketitle
\thispagestyle{empty}
{ \small \begin{center}
{\bf Abstract}
\end{center}

A family of $q$-difference-differential equations in two complex variables is studied, under the action of a so-called Mahler transform on time variable. The appearance of a leading formal $q$-difference operator of irregular type in the equation guarantees the existence of a formal solution to the main problem in time variable, which turns out to be obtained after a $q$-analog of Borel-Laplace procedure. Such formal solution to the main problem is $G_q$-summable along certain directions. The $G_q$-sums of such formal solutions do not in general satisfy the initial problem but rather turns out to be the analytic solutions to some related pseudo $q$-difference-differential equation. 

\smallskip

}
\bigskip \bigskip

\section{Introduction}

The present work studies a family of $q$-difference-differential equations in two variables in the complex domain. The equations in such family combine $q$-difference operators, Mahler transforms in time variable $t$, together with differential operators in $z$ variable. The precise shape of such family of equations are as follows: 
\begin{equation}\label{epralintro}
Q(\partial_z)u(t,z)=\exp_q\left(\alpha_Dt^{d_{D}}\sigma_{q;t}^{\frac{d_D}{k}}\right)R_D(\partial_z)u(t,z)+P(t,z,\sigma_{q;t},\partial_z,\{\mathfrak{m}_{\ell_2,t}\}_{\ell_2\in I})u(t,z)+f(t,z),
\end{equation}
under null initial data $u(0,z)\equiv 0$. In the previous equation, $q>1$ stands for a real number and $\sigma_{q,t}(f)(t)=f(qt)$ stands for the $q$-difference operator on $t$ variable, being the previous definition extended to $\sigma_{q,t}^{\delta}(f)(t)=f(q^{\delta}t)$ for any $\delta>0$. The elements $\alpha_D,d_D>0$ are positive real numbers and $k\ge1$ is an integer number. The terms $Q(X), R_D(X)$ belong to $\C[X]$ and $P(t,z,V_1,V_2,\{W_{\ell_2}\}_{\ell_2\in I})$ is some polynomial in the variables $t,V_1,V_2$ and a linear map in each $W_{\ell_2}$ for $\ell_2\in I$, for some fixed finite set $I\subseteq \N:=\{1,2,\ldots\}$. The function $P$ is holomorphic with respect to $z$ on a horizontal strip $H_{\beta}=\{z\in\C: |\hbox{Im}(z)|<\beta\}$, for some $\beta>0$.

The forcing term $f(t,z)$ is a polynomial in $t$ with bounded holomorphic coefficients defined in $H_{\beta}$. The symbol $\mathfrak{m}_{\ell_2,t}$, standing for the Mahler transform on $t$ variable, is defined by $\mathfrak{m}_{\ell_2,t}u(t,z)=u(t^{\ell_2},z)$ for $\ell_2\in I$.

Observe the appearance of a leading term in the equation consisting of a formal $q$-difference operator of infinite order. More precisely, the leading term involved an operator of the form
$$\exp_q(\alpha_Dt^{d_D}\sigma_{q;t}^{\frac{d_D}{k}})R_{D}(im)=R_D(im)\sum_{n\ge0}\frac{1}{[n]_q^{!}}\alpha_D^{n}(t^{d_D}\sigma_{q;t}^{\frac{d_D}{k}})^{(n)},$$
where $(t^{d_D}\sigma_{q;t}^{\frac{d_D}{k}})^{(n)}$ stands for the $n$-th iteration of the $q$-difference operator $t^{d_D}\sigma_{q;t}^{\frac{d_D}{k}}$.

\vspace{0.3cm}

This study puts forward the $q$-analog problem studied in~\cite{ma24}, where the second author constructs the formal solution to the nonlinear initial value problem with an infinite order irregular singularity and Mahler transform of the form
\begin{multline}\label{eauxintro2}
Q(\partial_z)u(t,z)=\hbox{cosh}\left(\alpha_D(t^{k+1}\partial_t)^2\right)R_D(\partial_z)u(t,z)+P(t,z,t^{k+1}\partial_t,\partial_z,\{\mathfrak{m}_{\ell_2,t}\}_{\ell_2\in I})u(t,z)\\
Q_1(\partial_z)u(t,z)\times Q_2(\partial_z)u(t,z)+f(t,z),
\end{multline} 
with $u(0,z)\equiv 0$. In the previous equation, $Q_1(X),Q_2(X)$ belong to $\C[X]$, with the rest of the elements being of the same nature as before.

\vspace{0.3cm}

The motivation of the problem considered in this paper is the previous research~\cite{ma24} arising from the interest and active research on the so-called Mahler equations. These are linear equations of the form
\begin{equation}\label{e01}
\sum_{p=0}^{n}b_p(z)v(z^{\ell^p})=0,
\end{equation}
for some integers $\ell\ge 2$, $n\ge1$ and rational coefficients $b_p(z)\in\C(z)$. We mention the applications of these equations to automatic sequences and Cobham's theorem (see~\cite{ada}),  Galoisian aspects and hypertranscendence (see~\cite{ada2}) or the algebraic structure of the solutions to the equation (\ref{e01}) (see~\cite{chy,fav,roq,roq2}). For a recent extension to the nonlinear setting, we refer to~\cite{bez,gongor25,nis}.

Our main equation not only involves Mahler transforms, but also formal operators of infinite order as its leading term. We provide some references on differential equations involving infinite order operators such as~\cite{ao,kast1,kast2,kast3}. In the framework of PDEs and hamiltonian systems related to physical problems, we refer to~\cite{deg,du,dub,dub2}.

In the work~\cite{scsi} there are different kind of coupling of operators acting on the equations of different nature: differential operators, $q$-difference, Mahler operator and shift operator. The general form of meromorphic solutions on the universal covering $\widetilde{\C\setminus\{0\}}$ are described. An extension of this result has been obtained in~\cite{ada2} by using the recent Galois theory developed by the last author of that paper.

\vspace{0.3cm}

In the previous work~\cite{ma24}, a formal solution $\hat{u}(t,z)$ to (\ref{eauxintro2}) as a formal power series in $t$ with coefficients being holomorphic functions in $H_{\beta}$ is constructed. Such formal power series is conformed by means of the classical Borel-Laplace method (see~\cite{ba} for a reference on this method). Such method associates $\hat{u}(t,z)$ to an analytic function $u^{d}(t,z)$ defined on $S \times H_{\beta}$, for some bounded sector $S$ with vertex at the origin and well-chosen bisecting direction. The fact that the opening of $S$ can be chosen to be wide enough guarantees unicity of the sum $u^{d}(t,z)$ having $\hat{u}(t,z)$ as its $k$-sum in $S$, with respect to $t$ variable. However, it is proved that the $k$-sum of the formal solution to (\ref{eauxintro2}) is not (in general) an analytic solution to (\ref{eauxintro2}) but different functional equations related to the initial problem.

The procedure considered in this work follows similar steps as in this previous research, with important variations arising from the $q$-nature of the equation under study, that need a different approach to be surpassed. We proceed to say some words on this concern.

In the search of formal solutions to equation (\ref{epralintro}), the strategy is to search for them in the form of an inverse Fourier transform followed by a second step based on the application of a formal $q$-Borel transformation of order $k$ to the equation (see Section~\ref{secI}). This procedure leads to an auxiliary functional equation in which the action of Mahler operator $\mathfrak{m}_{\ell_2,t}$ has turned into the following $q$-deceleration operator, constructed in Section~\ref{secb}:
\begin{equation}\label{e01}
\hat{\mathcal{D}}_{\ell_2}(\omega^d)(\tau^{\ell_2},m):=-\frac{q^{\frac{1}{8k'}}\sqrt{k'}}{(2\pi \log(q))^{1/2}}\int_{\tilde{C}_{R/2}}q^{\frac{k'}{2}\frac{\log(x/h)}{\log(q)}\left(\frac{\log(x/h)}{\log(q)}-\frac{1}{k'}\right)}\omega^d(\frac{x}{q^{k''}},m)\frac{dx}{x},
\end{equation}
where $\tilde{C}_{R/2}=\{R/2e^{i\theta}:\theta\in \R\}\subseteq\tilde{D}(0,R)$, for small enough $R>0$, and where $k'=\frac{k}{\ell_2^2-1}$ and $k''=\frac{\ell_2^2-\ell_2}{2k}$.

The strategy is based on the existence of a fixed point of certain operators in Banach spaces of functions whose definition and properties are described in Section~\ref{secbanach}. First, the fixed point is obtained in the Banach space of holomorphic functions near the origin with values in another Banach space related to the Fourier variable (Section~\ref{sec22}). The procedure is completed when considering the space of functions in the previous Banach space which can be analytically prolonged along an infinite sector with $q$-Gevrey bounds of certain order(Section~\ref{sec23}). It is important to proceed searching for the mentioned fixed point by working jointly on some neighborhood of the origin together with an infinite sector, due to the asymmetric behavior of the functions $\exp_q(z)$ and $\exp(\frac{\log^2|z|}{2\log(q)})$ near the origin and infinity. This last term appears in the operator (\ref{e01}). This is another difference observed regarding the procedure followed in the present work and~\cite{ma24}, where a symmetric behavior is observed. At this point, it is worth mentioning that the term $\hbox{cosh}(\cdot)$ at the head of equation (\ref{eauxintro2}) has the same behavior near infinity in the set of complex numbers of positive/negative real part, tending to infinity, which does not hold for the exponential function which has an asymmetric behavior there. On its side, the function $\exp_q(\cdot)$ approaches infinity in a domain which comprises the whole complex plane except from an infinite sector bisected by direction $\pi$ and of any small positive opening. This is a crucial property in the construction of the fixed point functions mentioned above.

Some auxiliary estimates on the geometric configuration of the problem is also described in Section~\ref{secgeo}, allowing to obtain the fixed point in the Banach spaces mentioned (Section~\ref{secbsi} and Section~\ref{secbsii}). The main results of the work are stated in Section~\ref{secmainres}. The first main result (Theorem~\ref{teo1}) states the existence of a formal solution $\hat{u}(t,z)$ to the main problem (\ref{epralintro}) in the form of a power series in $t$ with holomorphic and bounded functions on $H_{\beta'}$ for any fixed $0<\beta'<\beta$, which is $G_q$-summable of order $k$ with respect to some well-chosen direction $d\in\R$. See~\cite{zh99,mazh} for the notion of the concept of a $G_q$-summable formal power series, which in turn is a $q$-analog of classical $k$-summability as defined in~\cite{ba2,loday}. The second main result states that the $G_q$-sum of $\hat{u}(t,z)$, say $u^{d}(t,z)$, is holomorphic on $\tilde{D}\times H_{\beta'}$, where  $\tilde{D}$ stands for the universal covering of some neighborhood of the origin in the Riemann surface of the logarithm, $\mathcal{R}$. In addition to this, $u^{d}(t,z)$ satisfies a pseudo $q$-difference-differential and does not, in general, define an actual solution of the initial main problem under study (\ref{epralintro}), as it is justified in the last remark of this paper. The paper concludes with an appendix (Section~\ref{secapqex}) recalling known facts about $q$-exponential function, analytic $q$-Laplace and $q$-Borel transforms and their main properties, and a brief review on $G_q$-summability of formal power series.

\vspace{0.6cm}

\textbf{Notation:}

For $z_0\in\C$ and $r>0$, we write $D(z_0,r)$ for the open disc $\{z\in\C:|z-z_0|<r\}$, and $\overline{D}(z_0,r)$ for its closure. We write $\tilde{D}(0,r)$ for the universal covering of the disc $D(0,r)$ in the Riemann surface of the logarithm, $\mathcal{R}$.  

Let $\mathbb{E}$ be a Banach space of functions. Given an open set $U\subseteq\R$, $\mathcal{C}(U,\mathbb{E})$ stands for the set of continuous functions defined on $U$ with values in $\mathbb{E}$. We write $\mathcal{C}(U):=\mathcal{C}(U,\C)$. We write $\mathbb{E}\{\epsilon\}$ for the set of holomorphic functions defined on some neighborhood of the origin (in the variable $\epsilon$) and values on $\mathbb{E}$, and $\mathbb{E}[[\epsilon]]$ stands for the set of formal power series in $\epsilon$ with coefficients in $\mathbb{E}$. We denote by $\mathcal{O}_b(V,\mathbb{E})$ the set of bounded holomorphic functions on $V$  with values in $\mathbb{E}$. For simplicity, we write $\mathcal{O}(U)$ and $\mathcal{O}_b(V)$ instead of $\mathcal{O}(U,\C)$ and $\mathcal{O}_b(V,\C)$, respectively.

\section{Banach spaces of functions}\label{secbanach}

In this section, we recall the definition and main properties associated to certain Banach spaces used in the present work. In a first section, we deal with inverse Fourier transform, followed by a subsection considering functions defined on some neighborhood of the origin. In the third subsection, we consider a Banach space of functions which are holomorphic on a sectorial region  near infinity and subject to certain $q$-exponential growth.

\subsection{On inverse Fourier transform}\label{ap2}

Let $\beta>0$ and $\mu>1$. For any given $f\in\mathcal{C}(\R)$ satisfying there exists $C>0$ with 
\begin{equation}\label{e437}
|f(m)|\le C\frac{1}{(1+|m|)^{\mu}}\exp(-\beta|m|),\quad m\in\R,
\end{equation}
the inverse Fourier transform of $f$ is defined by
$$\mathcal{F}^{-1}(f)(x)=\frac{1}{(2\pi)^{1/2}}\int_{-\infty}^{\infty}f(m)\exp(imx)dm,\quad x\in\R.$$
The function $\mathcal{F}^{-1}(f)$ turns out to be holomorphic in the infinite horizontal strip 
$$H_{\beta}=\{z\in\C:|\hbox{Im}(z)|<\beta\},$$
and motivates the definition of the following Banach space.

\begin{defin}
For any continuous function $f:\R\to\C$ such that (\ref{e437}) holds, we define
$$\left\|h(m)\right\|_{(\beta,\mu)}:=\sup_{m\in\R}(1+|m|)^{\mu}\exp(\beta|m|)|f(m)|,$$
defining a norm, in the set of continuous functions satisfying (\ref{e437}), denoted by $E_{(\beta,\mu)}$. The pair $(E_{(\beta,\mu)},\left\|\cdot\right\|_{(\beta,\mu)})$ is a Banach space.
\end{defin}

The next properties hold regarding inverse Fourier transform. Their proof can be found in detail in~\cite{lama}, Proposition 7.

\begin{prop}
Let $\beta>0$, $\mu>1$. Given $h\in E_{(\beta,\mu)}$ then, the following statements hold:
\begin{itemize}
\item Define $\R\ni m\mapsto \varphi(m)=im h(m)$. Then, $\varphi\in E_{(\beta,\mu-1)}$ and $\partial_z\mathcal{F}^{-1}(h)(z)=\mathcal{F}^{-1}(\varphi)(z)$, for $z\in H_{\beta}$.
\item Given $g\in E_{(\beta,\mu)}$, the convolution of $h$ and $g$, denoted by $\R\ni m\mapsto \psi(m)=\frac{1}{(2\pi)^{1/2}}(h\star g)(m)$, and defined by
$$(h\star g)(m)=\int_{-\infty}^{\infty}h(m-m_1)g(m_1)dm_1$$
turns out to be an element of $E_{(\beta,\mu)}$, and it holds that
$$\mathcal{F}^{-1}(h)(z)\mathcal{F}^{-1}(g)(z)=\mathcal{F}^{-1}(\psi)(z),\quad z\in H_{\beta}.$$
\end{itemize}
\end{prop}





\subsection{Banach space of functions defined on some neighborhood of the origin}\label{sec22}

In this section, we describe an $\ell_1$-norm in spaces of convergent functions near the origin. We choose $\beta>0$ and $\mu>1$.

\begin{defin}
Let $R>0$. We define the norm $\left\|\cdot\right\|_{(1,R)}$, acting on the space of formal power series with values in $E_{(\beta,\mu)}$ as follows. Let $f(\tau)=\sum_{p\ge1}h_p\tau^p\in \tau E_{(\beta,\mu)}[[\tau]]$, with $h_p\in E_{(\beta,\mu)}$. We define the norm 
$$\left\|f\right\|_{(1,R)}:=\sum_{p\ge1}\left\|h_p\right\|_{(\beta,\mu)}R^{p}.$$
Let $HE_{(\beta,\mu,R)}$ be the Banach space of formal power series $f(\tau)$ as before such that $\left\|f\right\|_{(1,R)}<+\infty$.  
\end{defin}





\begin{lemma}\label{lema0}
Let $R>0$. Given $f\in HE_{(\beta,\mu,R)}$. Then, $f\in\mathcal{O}_b(D(0,R),E_{(\beta,\mu)})$. We have that 
$$\sup_{|\tau|< R}\left|f(\tau)\right|\le \left\|f(\tau)\right\|_{(1,R)}.$$
\end{lemma}

\begin{lemma}\label{lema3}
Let $A\in E_{(\beta,\mu)}$, $B(m)\in\C[m]$ with $B(0)\neq 0$. We choose integer numbers $h_0,h_1\ge0$. We also take $\alpha>0$ such that $\alpha\ge\hbox{deg}(B)$. Let $L(\tau,m)$ be a continuous function in $\overline{D}(0,\rho)\times \R$, holomorphic with respect to its first variable on $D(0,\rho)$ and such that if we write 
$$L(\tau,m)=\sum_{p\ge1}f_p(m)\tau^p,\quad \tau\in D(0,\rho),$$
then for all $0<R_1<\rho$, there exists $C=C(R_1)>0$ such that $|f_p(m)|\le \frac{C}{R_1^{p}(1+|m|)^{\alpha}}$ for all $m\in\R$ and all $p\ge0$.

Let $W(\tau,m)\in HE_{(\beta,\mu,R)}$, where $0<R<R_1$. Then, there exists $\tilde{C}_{1}>0$ such that
$$\left\|L(\tau,m)\int_{-\infty}^{\infty}A(m-m_1)\tau^{h_0}\sigma_{q;\tau}^{-h_1}W(\tau,m)B(im_1)dm_1\right\|_{(1,R)}\le \tilde{C}_1\left\| W(\tau,m)\right\|_{(1,R)}.$$  
\end{lemma}
\begin{proof}
We write $W(\tau,m)=\sum_{p\ge1}W_{p}(m)\tau^p$, for certain $W_p\in E_{(\beta,\mu)}$. It holds that
\begin{multline*}
\left\|L(\tau,m)\int_{-\infty}^{\infty}A(m-m_1)\tau^{h_0}\sigma_{q;\tau}^{-h_1}W(\tau,m)B(im_1)dm_1\right\|_{(1,R)}\\
=\sum_{p\ge1}\left\|\sum_{j=0}^{p}f_j(m)\int_{-\infty}^{\infty}A(m-m_1)\frac{W_{p-j}(m_1)}{q^{(p-j)h_1}}B(im_1)dm_1\right\|_{(\beta,\mu)}R^{p+h_0}\\
\le\sum_{p\ge1}\sum_{j=0}^{p}\left\|f_j(m)\int_{-\infty}^{\infty}A(m-m_1)\frac{W_{p-j}(m_1)}{q^{(p-j)h_1}}B(im_1)dm_1\right\|_{(\beta,\mu)}R^{p+h_0}
\end{multline*}

We observe there exists $C_{B}>0$ such that $|B(im)|\le C_B(1+|m|)^{\hbox{deg}(B)}$, valid for all $m\in\R$. In addition to this, taking into account the definition of $E_{(\beta,\mu)}$ and the hypothesis on the coefficients of Taylor expansion of $\tau\mapsto L(\tau,m)$, we arrive at

\begin{multline*}
\left\|f_j(m)\int_{-\infty}^{\infty}A(m-m_1)\frac{W_{p-j}(m_1)}{q^{(p-j)h_1}}B(im_1)dm_1\right\|_{(\beta,\mu)}\\
=\sup_{m\in\R}\left|f_j(m)\int_{-\infty}^{\infty}A(m-m_1)\frac{W_{p-j}(m_1)}{q^{(p-j)h_1}}B(im_1)dm_1\right|(1+|m|)^{\mu}e^{\beta|m|}\\ 
\le C C_B \sup_{m\in\R}\frac{1}{R_1^{j}(1+|m|)^{\alpha}}\int_{-\infty}^{\infty}|A(m-m_1)|\frac{|W_{p-j}(m_1)|}{q^{(p-j)h_1}}(1+|m_1|)^{\hbox{deg}(B)}dm_1(1+|m|)^{\mu}e^{\beta|m|}\\
\le C C_B \left\|A(m)\right\|_{(\beta,\mu)}\frac{1}{R_1^{j}}\frac{1}{q^{(p-j)h_1}}\left\|W_{p-j}(m)\right\|_{(\beta,\mu)}\\
\times\sup_{m\in\R}(1+|m|)^{\mu-\alpha}\int_{-\infty}^{\infty}\frac{1}{(1+|m-m_1|)^{\mu}}\frac{1}{(1+|m_1|)^{\mu-\hbox{deg}(B)}}dm_1.
\end{multline*}
In the previous upper estimates, we have made use of the inequality $e^{\beta|m|}e^{-\beta|m-m_1|}e^{-\beta|m_1|}\le 1$, coming from the triangular inequality.

We recall from Lemma 2.2~\cite{cota2} and the assumption that $\alpha>\hbox{deg}(B)$ that one can guarantee the existence of $C_{A,B}>0$ such that
$$\sup_{m\in\R}(1+|m|)^{\mu-\alpha}\int_{-\infty}^{\infty}\frac{dm_1}{(1+|m-m_1|)^{\mu}(1+|m_1|)^{\mu-\hbox{deg}(B)}}\le C_{A,B}.$$
This entails that
\begin{multline*}
\left\|f_j(m)\int_{-\infty}^{\infty}A(m-m_1)\frac{W_{p-j}(m_1)}{q^{(p-j)h_1}}B(im_1)dm_1\right\|_{(\beta,\mu)}\\
\le C C_{A,B} C_B \left\|A(m)\right\|_{(\beta,\mu)}\frac{1}{R_1^{j}}\frac{1}{q^{(p-j)h_1}}\left\|W_{p-j}(m)\right\|_{(\beta,\mu)}.
\end{multline*}
Therefore, one has that

\begin{multline*}
\left\|L(\tau,m)\int_{-\infty}^{\infty}A(m-m_1)\tau^{h_0}\sigma_{q;\tau}^{-h_1}W(\tau,m)B(im_1)dm_1\right\|_{(1,R)}\\
\le C C_{A,B} C_B \left\|A(m)\right\|_{(\beta,\mu)}\sum_{p\ge1}\sum_{j=0}^{p} \frac{1}{R_1^{j}}\frac{\left\|W_{p-j}(m_1)\right\|_{(\beta,\mu)}}{q^{(p-j)h_1}}R^{p+h_0}\\
= C C_{A,B} C_B R^{h_0} \left\|A(m)\right\|_{(\beta,\mu)}\left(\sum_{p\ge0}\left(\frac{R}{R_1}\right)^{p}\right)\left(\sum_{p\ge1}\left\|W_{p}(m)\right\|_{(\beta,\mu)}\frac{R^{p}}{q^{ph_1}}\right)\\
\le C C_{A,B} C_B R^{h_0} \left\|A(m)\right\|_{(\beta,\mu)}\left(\sum_{p\ge0}\left(\frac{R}{R_1}\right)^{p}\right)\left\|W(\tau,m)\right\|_{(1,R)}.
\end{multline*}

The result follows by taking 
$$\tilde{C}_1=C C_{A,B} C_B R^{h_0} \left\|A(m)\right\|_{(\beta,\mu)}\left(\sum_{p\ge0}\left(\frac{R}{R_1}\right)^{p}\right).$$
\end{proof}

In the following result, we describe the action of certain operator involving a $q$-deceleration operator on the previous Banach space.  

\begin{defin}\label{def216}
Let $\hat{f}(\tau)=\sum_{n\ge1}f_n\tau^n\in \tau E_{(\beta,\mu)}[[\tau]]$. Given an integer $p\ge2$, we define the formal $q$-deceleration operator by
$$\hat{\mathcal{D}}_p(\hat{f})(h)=\sum_{n\ge1}f_n\frac{q^{\frac{n(n-1)}{2k}}}{q^{\frac{pn(pn-1)}{2k}}}h^n.$$
\end{defin}

\begin{lemma}\label{lema4}
Let $A\in E_{(\beta,\mu)}$, $B(m)\in\C[m]$ with $B(0)\neq 0$. We also take $\alpha\ge0$ such that $\alpha\ge\hbox{deg}(B)$, and choose integer numbers $h_0\ge0$ ,$h_1\ge 1$ and $h_2\ge2$. Let $L(\tau,m)$ be a continuous function in $\overline{D}(0,\rho)\times \R$, holomorphic with respect to its first variable on $D(0,\rho)$ and such that if we write
$$L(\tau,m)=\sum_{p\ge0}f_p(m)\tau^p,\quad \tau\in D(0,\rho),$$
then for all $0<R_1<\rho$, there exists $C=C(R_1)>0$ such that $|f_p(m)|\le \frac{C}{R_1^{p}(1+|m|)^{\alpha}}$ for all $m\in\R$ and all $p\ge0$.

Let $W(\tau,m)\in HE_{(\beta,\mu,R)}$, and $0<R<R_1$. 

Then, for small enough $0<R<R_1$, there exists $\tilde{C}_{2}>0$ such that
$$\left\|L(\tau,m)\int_{-\infty}^{\infty}A(m-m_1)\hat{\mathcal{D}}_{h_2}\left(\xi\mapsto \xi^{h_0}\sigma_{q;\xi}^{-h_1}W(\xi,m_1)\right)(\tau^{h_2})B(im_1)dm_1\right\|_{(1,R)}\le \tilde{C}_2\left\| W(\tau,m)\right\|_{(1,R)}.$$  
\end{lemma}
\begin{proof}
Let $W(\tau,m)=\sum_{p\ge 1}W_{p}(m)\tau^p$, for certain $W_p\in E_{(\beta,\mu)}$. We observe that
$$
\hat{\mathcal{D}}_{h_2}\left(\xi\mapsto \xi^{h_0}\sigma_{q;\xi}^{-h_1}W(\xi,m_1)\right)(\tau^{h_2})
=\sum_{n\ge h_0}W_{n-h_0}(m_1)\frac{q^{\frac{n(n-1)}{2k}}}{q^{\frac{h_2n(h_2n-1)}{2k}}}q^{-h_1(n-h_0)}\tau^{nh_2}.$$
Therefore, one has
\begin{multline*}
\left\|L(\tau,m)\int_{-\infty}^{\infty}A(m-m_1)\hat{\mathcal{D}}_{h_2}\left(\xi\mapsto \xi^{h_0}\sigma_{q;\xi}^{-h_1}W(\xi,m)\right)(\tau^{h_2})B(im_1)dm_1\right\|_{(1,R)}\\
=\sum_{n\ge h_0}\left\|\sum_{n_1+n_2=n}f_{n_1}(m)\int_{-\infty}^{\infty}A(m-m_1)W_{n_2-h_0}(m_1)q^{-h_1(n_2-h_0)}\frac{q^{\frac{n_2(n_2-1)}{2k}}}{q^{\frac{h_2n_2(h_2n_2-1)}{2k}}}\tau^{n_1}\tau^{h_2n_2}B(im_1)dm_1\right\|_{(\beta,\mu)}\\
\le \sum_{n\ge h_0}\left(\sum_{n_1+n_2=n}\left\|f_{n_1}(m)\int_{-\infty}^{+\infty}A(m-m_1)W_{n_2-h_0}(m_1)q^{-h_1(n_2-h_0)} B(im_1)dm_1\right\|_{(\beta,\mu)}\right.\\
\left.\times\frac{q^{\frac{n_2(n_2-1)}{2k}}}{q^{\frac{n_2h_2(n_2h_2-1)}{2k}}}\right)R^{n_2+h_2n_2}\\
\le \sum_{n\ge h_0}\sum_{n_1+n_2=n}(1+|m|)^{\mu}e^{\beta|m|}\frac{C}{R_1^{n_1}(1+|m|)^{\alpha}}\int_{-\infty}^{+\infty}\left\|A(m)\right\|_{(\beta,\mu)}\frac{1}{(1+|m-m_1|)^{\mu}}\\
\times e^{-\beta|m-m_1|}\left\|W_{n_2-h_0}\right\|_{(\beta,\mu)}(1+|m_1|)^{-\mu}e^{-\beta|m_1|}q^{-h_1(n_2-h_0)}C_B(1+|m_1|)^{\hbox{deg}(B)}dm_1\\
\times  \frac{q^{\frac{n_2(n_2-1)}{2k}}}{q^{\frac{n_2h_2(n_2h_2-1)}{2k}}}R^{n_2+h_2n_2}\\
\le C C_{B}\left\|A(m)\right\|_{(\beta,\mu)}C_{A,B}\sum_{n\ge h_0}\sum_{n_1+n_2=n}\frac{1}{R^{n_1}}\left\|W_{n_2-h_0}(\tau,m)\right\|_{(\beta,\mu)}q^{-h_1(n_2-h_0)}\frac{q^{\frac{n_2(n_2-1)}{2k}}}{q^{\frac{n_2h_2(n_2h_2-1)}{2k}}}R^{n_2+h_2n_2}.
\end{multline*}
At this point, we observe there exists $C_q>0$ such that
$$\frac{q^{\frac{n_2(n_2-1)}{2k}}}{q^{\frac{n_2h_2(n_2h_2-1)}{2k}}}q^{-h_1(n_2-h_0)}=q^{\epsilon_1 n_2^2+\epsilon_2 n_2+\epsilon_3}\le C_q,$$
where $\epsilon_1=(1-h_2^2)/(2k)<0$, $\epsilon_2=(h_2-1)/(2k)-h_1$ and $\epsilon_3=h_1h_0$, valid  for all $n_2\ge0$. This yields
\begin{multline*}
\left\|L(\tau,m)\int_{-\infty}^{\infty}A(m-m_1)\hat{\mathcal{D}}_{h_2}\left(\xi\mapsto \xi^{h_0}\sigma_{q;\xi}^{-h_1}W(\xi,m)\right)(\tau^{h_2})B(im_1)dm_1\right\|_{(1,R)}\\
\le C C_{B}\left\|A(m)\right\|_{(\beta,\mu)}C_{A,B}C_q\sum_{n\ge h_0}\sum_{n_1+n_2=n}\frac{1}{R_1^{n_1}}\left\|W_{n_2-h_0}\right\|_{(\beta,\mu)}R^{n_2+h_2n_2}\\
\le C C_{B}\left\|A(m)\right\|_{(\beta,\mu)}C_{A,B}C_qR^{h_0}\sum_{n\ge h_0}\sum_{n_1+n_2=n}\frac{1}{R_1^{n_1}}\left\|W_{n_2-h_0}\right\|_{(\beta,\mu)}R^{n_1+n_2-h_0}\\
\le C C_{B}\left\|A(m)\right\|_{(\beta,\mu)}C_{A,B}C_qR^{h_0}\left(\sum_{n\ge0}\left(\frac{R}{R_1}\right)^{n}\right)\left(\sum_{n\ge0} \left\|W_{n}(m)\right\|_{(\beta,\mu)}R^n\right)
\end{multline*}
The result follows for
$$\tilde{C}_2=CC_BC_{A,B}C_q\left\|A(m)\right\|_{(\beta,\mu)}R^{h_0}\left(\sum_{p\ge 0}\left(\frac{R}{R_1}\right)^{p}\right).$$
\end{proof}

\label{remark}\textbf{Remark:} Observe that the constants $\tilde{C}_1,\tilde{C}_2$ in Lemma~\ref{lema3} and Lemma~\ref{lema4}, respectively, can be made arbitrary small when considering $R>0$ sufficiently small, with $h_0\ge 1$.

\subsection{Formal $q$-Borel transform of order $k$. Mahler transform and $q$-deceleration operator}\label{secb}

In this section, we define the formal $q$-Borel transform of order $k>0$, which is intimately related to its analytic counterpart analyzed in Section~\ref{sec332}. Later, we define Mahler transform and its influence on the action on the previous formal operator, which will be crucial in the next section in order to follow with a second step in the strategy to search for solutions to the main problem.

In the whole subsection, $(\mathbb{E},\left\|\cdot\right\|_{\mathbb{E}})$ stands for a complex Banach space, and $k\ge1$ is an integer.

Regarding Corollary~\ref{coro467} in the Appendix, we define the formal Borel transform of order $k$ as follows.

\begin{defin}
Let $\hat{U}(T)=\sum_{n\ge 1}a_nT^n\in T\mathbb{E}[[T]]$. The formal $q$-Borel transform of order $k$ is defined by
\begin{equation}\label{e713}
\hat{\mathcal{B}}_{q;k}(U(T))(\xi)=\sum_{n\ge1}\frac{a_n}{q^{\frac{n(n-1)}{2k}}}\xi^n.
\end{equation}
\end{defin}

The following is a straightforward result derived from the definition.

\begin{prop}\label{prop195}
Let $\sigma,j\ge0$ be integer numbers and let $\hat{U}(T)\in T\mathbb{E}[[T]]$. Then, one has that
$$\hat{\mathcal{B}}_{q;k}(T^{\sigma}\sigma_{q;T}^j(\hat{U}(T)))(\xi)=\frac{\xi^\sigma}{q^{\frac{\sigma(\sigma-1)}{2k}}}\sigma_{q;\xi}^{j-\frac{\sigma}{k}}\hat{\mathcal{B}}_{q;k}(\hat{U}(T))(\xi).$$
\end{prop}
\begin{proof}
Let $\hat{U}(T)=\sum_{n\ge1}a_nT^n$. By definition, one has that
$$\hat{\mathcal{B}}_{q;k}(T^{\sigma}\sigma_{q;T}^j(\hat{U}(T)))(\xi)=\sum_{n\ge1}a_nq^{jn}\frac{1}{q^{\frac{(n+\sigma)(n+\sigma-1)}{2k}}}\xi^{n+\sigma}.$$
On the other hand, 
$$\sigma_{q;\xi}^{j-\frac{\sigma}{k}}\hat{\mathcal{B}}_{q;k}(\hat{U}(T))(\xi)=\sum_{j\ge1}a_n\frac{1}{q^{\frac{n(n-1)}{2k}}}q^{(j-\frac{\sigma}{k})n}\xi^n,$$
which leads us to the conclusion, after observing that
$$\frac{(n+\sigma)(n+\sigma-1)}{2k}=\frac{\sigma(\sigma-1)}{2k}+\frac{n(n-1)}{2k}+\frac{2\sigma n}{2k}.$$
\end{proof}

At this point, we analyze the action of the so-called Mahler transform $T\mapsto T^p$ for some fixed $p\in\N$, $p\ge2$, on the formal Borel map.

\begin{lemma}
Given the formal power series $\hat{U}(T)=\sum_{n\ge1}a_nT^n\in T\mathbb{E}[[T]]$, we write $\hat{V}(T)=\hat{U}(T^p)$. Then, it holds that
$$\hat{\mathcal{B}}_{q;k}(\hat{V})(\xi)=\sum_{n\ge1}a_n\frac{\xi^{pn}}{q^{\frac{pn(pn-1)}{2k}}}.$$
\end{lemma}

The relationship between $\hat{\mathcal{B}}_{q;k}(\hat{V})(\xi)$ and $\hat{\mathcal{B}}_{q;k}(\hat{U})(\xi)$ is given in terms of a $q$-analog of the so-called deceleration operator defined in~\cite{ba}, page 46. See also~\cite{ma24}, which was stated in Definition~\ref{def216}.

\begin{prop}\label{prop232}
Given the formal power series $\hat{U}(T)=\sum_{n\ge1}a_nT^n\in T\mathbb{E}[[T]]$, and let $\hat{V}(T)=\hat{U}(T^p)$. Then, one has that
$$\hat{\mathcal{B}}_{q;k}(\hat{V}(T))(\xi)=\hat{\mathcal{D}}_{p}(\hat{\mathcal{B}}_{q;k}(\hat{U}))(\xi^p).$$
\end{prop}
\begin{proof}
It is straightforward from the definition of the formal $q$-Borel and formal $q$-deceleration operators acting on $\tau E[[\tau]]$, defined in an analogous way as in Definition~\ref{def216} for the specific Banach space $E_{(\beta,\mu)}$.
 \end{proof}

It is natural to extend Definition~\ref{def216} to analytic functions defined on some neighborhood the origin by identifying the function with their Taylor representation at the origin. Then, the following result holds. Now, we provide an integral representation for the formal $q$-deceleration operator under analytic assumptions on the formal power series $\hat{f}$ in Definition~\ref{def216}.

\begin{prop}\label{prop233}
Let $p\ge2$ be an integer number. We fix a holomorphic function $f$ defined on some neighborhood of the origin, and with values in $\mathbb{E}$. Then, one has:
\begin{itemize}
\item $\hat{\mathcal{D}}_{p}(f)$ defines an entire function.

Let us define $k'=\frac{k}{p^2-1}>0$ and $k''=\frac{p^2-p}{2k}>0$.

\item For all $h\in\C\setminus\{0\}$ one has
$$\hat{\mathcal{D}}_p(f)(h)=-\frac{q^{\frac{1}{8k'}}\sqrt{k'}i}{(2\pi \log(q))^{1/2}}\int_{\tilde{C}_{R/2}}q^{\frac{k'}{2}\frac{\log(x/h)}{\log(q)}\left(\frac{\log(x/h)}{\log(q)}-\frac{1}{k'}\right)}\hat{f}(\frac{x}{q^{k''}})\frac{dx}{x},$$
where $\tilde{C}_{R/2}=\{R/2e^{i\theta}:\theta\in \R\}\subseteq\tilde{D}(0,R)$, for small enough $R>0$.
\item There exist $C,K>0$ and $\Delta>1$ such that
\begin{equation}\label{e244}
\left|\hat{\mathcal{D}}_p(f)(h)\right|\le K\left(\sup_{x\in\tilde{C}_{R/2}}|\hat{f}(x/q'')|\right)\exp\left(\frac{k'}{2\log(q)}\log^2(|h|+\Delta)+\alpha\log(|h|+\Delta)\right),
\end{equation}
for all $h\in\C$ with $|h|\ge 1$.
\end{itemize}
\end{prop}
\begin{proof}
We make the assumption that $f$ is convergent on some disc $D(0,\rho)$ for some $\rho>0$. We rewrite
$$\frac{q^{\frac{n(n-1)}{2k}}}{q^{\frac{pn(pn-1)}{2k}}}=\frac{1}{q^{\frac{n(n-1)}{2k'}}}\frac{1}{q^{k''n}},$$
for $k'=\frac{k}{p^2-1}$ and $k''=\frac{p^2-p}{2k}>0$, if $p\ge 2$. From Corollary~\ref{coro467}, we can write
\begin{multline*}
\hat{\mathcal{D}}_p(f)(h)=\sum_{n\ge1}f_n\frac{1}{q^{k''n}}\frac{1}{q^{\frac{n(n-1)}{2k'}}}h^{n}=\sum_{n\ge1}f_n\frac{1}{q^{k''n}}\mathcal{B}_{q;k'}^{d}(\tau^n)(h)=\mathcal{B}_{q;k'}^{d}\left(\sum_{n\ge1}f_n\left(\frac{1}{q^{k''}}\tau\right)^n\right)(h)\\
=\mathcal{B}_{q;k'}^{d}(f(\tau/q^{k''}))(h).
\end{multline*}
At this point, we refer to Section~\ref{sec332}, on the properties of the analytic Borel transform.

The second part of the result is a direct consequence of Proposition~\ref{prop302}. Observe that the mapping $x\mapsto f(x/q^{k''})$ is univalued on $D(0,R)$. Observe moreover that the estimates in (\ref{e415}) are valid for all $t,\theta\in\R$, for well chosen $K>0$ according to the fact that $x\mapsto f(x/q^{k''})$ is univalued on a neighborhood of the origin.
\end{proof}

As a consequence of Proposition~\ref{prop232} and Proposition~\ref{prop233} we get the following result.

\begin{prop}\label{prop234}
Let $\hat{U}(T)\in T\mathbb{E}[[T]]$ and define $\hat{V}(T)=\hat{U}(T^p)$ for some integer $p\ge2$. Assume that $\mathcal{B}_{q;k}(\hat{U})(\tau)$ is convergent on some disc $D(0,\rho)$, for some $\rho>0$. Then, $\mathcal{B}_{q;k}(\hat{V})(\tau)$ is analytic in $\C$ and satisfies there exist $K,\alpha>0$ and $\Delta>1$ such that
$$\left|\mathcal{B}_{q;k}(\hat{V})(\xi)\right|\le K\exp\left(\frac{k'}{2\log(q)}\log^2(|\xi|^p+\Delta)+\alpha\log(|\xi|^p+\Delta)\right),$$
for all $\xi\in \C$, $|\xi|\ge 1$, where $k'=k/(p^2-1)$, or equivalently, there exist $\tilde{K},\tilde{\alpha}>0$ and $\Delta>1$
$$\left|\mathcal{B}_{q;k}(\hat{V})(\xi)\right|\le \tilde{K}\exp\left(\frac{k}{2\log(q)}\left(\frac{p^2}{p^2-1}\right)\log^2(|\xi|+\Delta)+\tilde{\alpha}\log(|\xi|+\Delta)\right),$$
for all $\xi\in\C$ with $|\xi|\ge 1$.
\end{prop}
\begin{proof}
The result is a direct consequence of Proposition~\ref{prop232} and Proposition~\ref{prop233}. The last equivalent result comes from the fact that  for $|\xi|\ge 1$ one has that
$$\log(|\xi|^p+\Delta)=\log(|\xi|^p(1+\frac{\Delta}{|\xi|^p}))=\log(|\xi|^p)+\log(1+\frac{\Delta}{|\xi|^p})\le p\log|\xi|+C\le p\log(|\xi|+\Delta)+C,$$
for some $C>0$.
\end{proof}

\subsection{Banach space of functions of $q$-exponential growth at infinity}\label{sec23}

We consider in this section the following Banach space of functions, as a modification of that appearing in~\cite{malek20}, based on the functional spaces developed in~\cite{drelasmal,lama12,lama15}. In the whole section, we fix $q>1$.

Given $\alpha,\rho>0$ we define the following Banach space of functions holomorphic on some neighborhood of the origin, and of $q$-exponential growth at infinity in the infinite sector $S_d$.

\begin{defin}\label{defi111}
Let $R>0$. The set $E_{(R,k,\alpha,\beta,\mu)}^d$ consists of all formal power series in $HE_{(\beta,\mu,R)}$ which satisfy the following property: from Lemma~\ref{lema0} such a formal power series defines a bounded holomorphic function on $D(0,R)$ with values in $E_{(\beta,\mu)}$, say $h(\tau,m)$. In addition to this, we assume such function $h(\tau,m)$ can be extended as a continuous function on $(S_d\cup D(0,R))\times\R$, holomorphic with respect to $\tau$ on $S_d\cup D(0,R)$ such that
$$\left\|h(\tau,m)\right\|_{(R,k,\alpha,\beta,\mu)}:=\left\|h(\tau,m)\right\|_{(1,R)}+\left\|h(\tau,m)\right\|_{(k,\alpha,\beta,\mu)}$$
is finite, with
\begin{multline*}
\left\|h(\tau,m)\right\|_{(k,\alpha,\beta,\mu)}:= \sup_{(\tau,m)\in  S_d\times \R}(1+|m|)^{\mu}e^{\beta|m|}\frac{1}{|\tau|}\\
\times\exp\left(-\frac{k}{2\log(q)}\log^2|\tau|-\alpha \log|\tau|\right)|h(\tau,m)|.
\end{multline*}
The set $(E^d_{(R,k,\alpha,\beta,\mu)},\left\|\cdot\right\|_{(R,k,\alpha,\beta,\mu)})$ is a Banach space over $\C$. 
\end{defin}

\textbf{Remark:} It is worth remarking that the formal power series in the previous Banach space determine functions which are of $q$-exponential growth near infinity with respect to variable $\tau$. However, they remain bounded near the origin, due to the finiteness of $\left\|\cdot\right\|_{(1,R)}$.

\begin{lemma}\label{lema486a}
Let $A\in E_{(\beta,\mu)}$, $B(m)\in\C[m]$ with $B(0)\neq 0$.  Let $d\in\R$. We choose integer numbers $h_0,h_1\ge 1$ and $0<R<\rho$. We also take $\tilde{\alpha}>0$ such that $\tilde{\alpha}\ge \hbox{deg}(B)$. Let $L(\tau,m)$ be a continuous function in $\overline{D}(0,\rho)\times \R$, holomorphic with respect to its first variable on $D(0,\rho)$ and such that if we write
$$L(\tau,m)=\sum_{p\ge0}f_p(m)\tau^p,\quad \tau\in D(0,\rho),$$
then for all $0<R_1<R$, there exists $C_3=C(R_1)>0$ such that $|f_p(m)|\le \frac{C_3}{R_1^{p}(1+|m|)^{\alpha}}$ for all $m\in\R$ and all $p\ge0$. Moreover, we assume this function can be analytically extended to an infinite sector of bisecting direction $d$, say $S_{d}$, on its first variable, satisfying that there exist $C_1,\tilde{K}>0$ and $h\in\R$ with
$$|L(\tau,m)|\le C_1\exp\left(-\frac{\tilde{K}\log^2|\tau|}{2\log(q)}\right)|\tau|^{h}\frac{1}{(1+|m|)^{\tilde{\alpha}}},$$
for all $(\tau,m)\in \overline{S_{d}}\times \R$.

Then, for every $0<R<\rho$, there exists $\tilde{C}_3>0$ such that for all $W(\tau,m)\in E^d_{(R,k,\alpha,\beta,\mu)}$, one has that
\begin{equation}\label{e343}
\left\|L(\tau,m)\int_{-\infty}^{\infty}A(m-m_1)\tau^{h_0}\sigma_{q;\tau}^{-h_1}W(\tau,m_1)B(im_1)dm_1\right\|_{(R,k,\alpha,\beta,\mu)}\le\tilde{C}_3 \left\| W(\tau,m)\right\|_{(R,k,\alpha,\beta,\mu)}.
\end{equation} 
\end{lemma}

\begin{proof}
Let $W(\tau,m)\in E^d_{(R,k,\alpha,\beta,\mu)}$. We divide the proof in two steps:

Step 1: 
The proof of Lemma~\ref{lema3} can be rewritten to arrive at 
\begin{equation}\label{e344}
\left\|L(\tau,m)\int_{-\infty}^{\infty}A(m-m_1)\tau^{h_0}\sigma_{q;\tau}^{-h_1}W(\tau,m_1)B(im_1)dm_1\right\|_{(1,R)}\le\tilde{C}_1\left\|W(\tau,m)\right\|_{(1,R)}.
\end{equation}

Step 2: Assume that $\tau\in S_d$. 

Step 2.1: Assume that $\tau\in S_d$ with $|\tau|<R$. Then, an analogous reasoning as that in Step 1, and from the fact that $\exp\left(-\frac{k}{2\log(q)}\log^2|\tau|-\alpha\log|\tau|\right)$ is upper bounded, and $|\tau|^{h_0-1}\le R_1^{h_0-1}$,  yield
\begin{multline}\label{e355}
\left|L(\tau,m)\int_{-\infty}^{\infty}A(m-m_1)\tau^{h_0}\sigma_{q;\tau}^{-h_1}W(\tau,m_1)B(im_1)dm_1\right|(1+|m|)^{\mu}e^{\beta|m|}\frac{1}{|\tau|}\\
\times\exp\left(-\frac{k}{2\log(q)}\log^2|\tau|-\alpha\log|\tau|\right)\\
\le C_B C_3\left(\sum_{p\ge0}\left(\frac{R}{R_1}\right)^p\right)C_{A,B}\left\|A(m)\right\|_{(\beta,\mu)}R^{h_0-1}\left\|W(\tau,m)\right\|_{(1,R)}.
\end{multline}

Step 2.2. Assume that $\tau\in S_d$ with $|\tau|\ge R$. Then, one has that
\begin{multline*}
\left|L(\tau,m)\int_{-\infty}^{\infty}A(m-m_1)\tau^{h_0}\sigma_{q;\tau}^{-h_1}W(\tau,m_1)B(im_1)dm_1\right|(1+|m|)^{\mu}e^{\beta|m|}\frac{1}{|\tau|}\\
\times\exp\left(-\frac{k}{2\log(q)}\log^2|\tau|-\alpha\log|\tau|\right)\\
\le \frac{C_1}{(1+|m|)^{\tilde{\alpha}}}\frac{|\tau|^{h}\exp\left(-\frac{k}{2\log(q)}\log^2|\tau|-\alpha\log|\tau|\right)}{\exp\left(\frac{\tilde{K}\log^2|\tau|}{2\log(q)}\right)}(1+|m|)^{\mu}e^{\beta|m|}\frac{1}{|\tau|}\\
\times\int_{-\infty}^{+\infty}\left\|A(m)\right\|_{(\beta,\mu)}\frac{1}{(1+|m-m_1|)^{\mu}}e^{-\beta|m-m_1|}|\tau|^{h_0+h}\left\|W(\tau,m)\right\|_{(k,\alpha,\beta,\mu)}\frac{1}{(1+|m_1|)^{\mu}}e^{-\beta|m_1|}\left|\frac{\tau}{q^{h_1}}\right|\\
\times \exp\left(\frac{k}{2\log(q)}\log^2\left|\frac{\tau}{q^{h_1}}\right|+\alpha\log\left|\frac{\tau}{q^{h_1}}\right|\right)C_B(1+|m_1|)^{\hbox{deg}(B)}dm_1\\
\le C_1 C_{A,B}\left\|A(m)\right\|_{(\beta,\mu)}\frac{1}{q^{h_1}}\left\|W(\tau,m)\right\|_{(k,\alpha,\beta,\mu)}\mathcal{A}_1(|\tau|),
\end{multline*}
with 
$$\mathcal{A}_1(|\tau|)=|\tau|^{h_0+h}\frac{\exp\left(-\frac{k}{2\log(q)}\log^2|\tau|-\alpha\log|\tau|\right)}{\exp\left(\frac{\tilde{K}\log^2|\tau|}{2\log(q)}\right)}\exp\left(\frac{k}{2\log(q)}\log^2\left|\frac{\tau}{q^{h_1}}\right|+\alpha\log\left|\frac{\tau}{q^{h_1}}\right|\right).$$
It is straight to check that $\sup_{|\tau|\ge R}\mathcal{A}_1(|\tau|)\le C_{\mathcal{A}_1}$, for some $C_{\mathcal{A}_1}>0$. This entails that
\begin{multline}\label{e373}
\left|L(\tau,m)\int_{-\infty}^{\infty}A(m-m_1)\tau^{h_0}\sigma_{q;\tau}^{-h_1}W(\tau,m_1)B(im_1)dm_1\right|(1+|m|)^{\mu}e^{\beta|m|}\frac{1}{|\tau|}\\
\times\exp\left(-\frac{k}{2\log(q)}\log^2|\tau|-\alpha\log|\tau|\right)\\
\le C_{\mathcal{A}_1}C_1 C_{A,B}\left\|A(m)\right\|_{(\beta,\mu)}\frac{1}{q^{h_1}}\left\|W(\tau,m)\right\|_{(k,\alpha,\beta,\mu)}.
\end{multline}

In view of (\ref{e355}) and (\ref{e373}) we deduce that 
\begin{multline}\label{e380}
\left\|L(\tau,m)\int_{-\infty}^{\infty}A(m-m_1)\tau^{h_0}\sigma_{q;\tau}^{-h_1}W(\tau,m_1)B(im_1)dm_1\right\|_{(k,\alpha,\beta,\mu)}\\
\le \left\|A(m)\right\|_{(\beta,\mu)}C_{A,B}\max\left\{C_B C_3\left(\sum_{p\ge0}\left(\frac{R}{R_1}\right)^p\right)R^{h_0-1}\left\|W(\tau,m)\right\|_{(1,R)},C_{\mathcal{A}_1}C_1\frac{1}{q^{h_1}}\left\|W(\tau,m)\right\|_{(k,\alpha,\beta,\mu)}\right\}.
\end{multline}

Finally, taking into account (\ref{e344}) and (\ref{e380}) one achieves (\ref{e343}) with
$$\tilde{C}_3=\max\{\tilde{C}_1+C_B C_3\left(\sum_{p\ge0}\left(\frac{R}{R_1}\right)^p\right)C_{A,B}\left\|A(m)\right\|_{(\beta,\mu)}R^{h_0-1},\tilde{C}_1+C_{\mathcal{A}_1}C_1 C_{A,B}\left\|A(m)\right\|_{(\beta,\mu)}\frac{1}{q^{h_1}}\}.$$
\end{proof}

\begin{lemma}\label{lema486b}
Let $A\in E_{(\beta,\mu)}$, $B(m)\in\C[m]$ with $B(0)\neq 0$.  Let $d\in\R$. We choose integer numbers $h_0\ge 1$, $h_1\ge 1$, $h_2\ge 2$, and $0<R<\rho$. We also take $\tilde{\alpha}>0$ such that $\tilde{\alpha}\ge \hbox{deg}(B)$.

Let $L(\tau,m)$ be a continuous function in $\overline{D}(0,\rho)\times \R$, holomorphic with respect to its first variable on $D(0,\rho)$ and such that if we write
$$L(\tau,m)=\sum_{p\ge0}f_p(m)\tau^p,\quad \tau\in D(0,\rho),$$
then for all $0<R_1<R$, there exists $C_3=C(R_1)>0$ such that $|f_p(m)|\le \frac{C_3}{R_1^{p}(1+|m|)^{\alpha}}$ for all $m\in\R$ and all $p\ge0$. Moreover, we assume this function can be analytically extended to an infinite sector of bisecting direction $d$, say $S_{d}$, on its first variable, satisfying that there exist $C_1,\tilde{K}>0$ and $h\in\R$ with
$$|L(\tau,m)|\le C_1\exp\left(-\frac{\tilde{K}\log^2|\tau|}{2\log(q)}\right)|\tau|^{h}\frac{1}{(1+|m|)^{\tilde{\alpha}}},$$
for all $(\tau,m)\in \overline{S_{d}}\times \R$, where
$$\tilde{K}>\frac{k}{h_2^2-1}.$$

Then, there exists $\tilde{C}_4>0$ such that for all $W(\tau,m)\in E^d_{(R,k,\alpha,\beta,\mu)}$, one has that
\begin{multline}\label{e343b}
\left\|L(\tau,m)\int_{-\infty}^{\infty}A(m-m_1)\hat{\mathcal{D}}_{h_2}\left(\xi\mapsto \xi^{h_0}\sigma_{q;\xi}^{-h_1}W(\xi,m)\right)(\tau^{h_2})B(im_1)dm_1\right\|_{(R,k,\alpha,\beta,\mu)}\\
\le\tilde{C}_4 \left\| W(\tau,m)\right\|_{(R,k,\alpha,\beta,\mu)}.
\end{multline}
\end{lemma}

\begin{proof}
Let $W(\tau,m)\in E^d_{(R,k,\alpha,\beta,\mu)}$. We divide the proof in two steps:

Step 1: 
The proof of Lemma~\ref{lema4} can be rewritten to arrive at 
\begin{equation}\label{e344b}
\left\|L(\tau,m)\int_{-\infty}^{\infty}A(m-m_1)\hat{\mathcal{D}}_{h_2}\left(\xi\mapsto \xi^{h_0}\sigma_{q;\xi}^{-h_1}W(\xi,m_1)\right)(\tau^{h_2})B(im_1)dm_1\right\|_{(1,R)}\le\tilde{C}_2\left\|W(\tau,m)\right\|_{(1,R)}.
\end{equation}

Step 2: Assume that $\tau\in S_d$. 

Step 2.1: Assume that $\tau\in S_d$ with $|\tau|<R$. Then, one considers
\begin{multline}\label{e421bb}
\left|L(\tau,m)\int_{-\infty}^{\infty}A(m-m_1)\hat{\mathcal{D}}_{h_2}\left(\xi\mapsto \xi^{h_0}\sigma_{q;\xi}^{-h_1}W(\xi,m_1)\right)(\tau^{h_2})B(im_1)dm_1\right|(1+|m|)^{\mu}e^{\beta|m|}\frac{1}{|\tau|}\\
\times \exp\left(-\frac{k}{2\log(q)}\log^2|\tau|-\alpha\log|\tau|\right).
\end{multline}
We observe that for $|\tau|<R$ one has
$$|L(\tau,m)|\le \sum_{p\ge0}\frac{C_3}{R_1^{p}(1+|m|)^{\tilde{\alpha}}}|\tau|^{p}\le C_3\sum_{p\ge0}\left(\frac{R}{R_1}\right)^{p}\frac{1}{(1+|m|)^{\tilde{\alpha}}}.$$
On the other hand,
$$\hat{\mathcal{D}}_{h_2}(\xi\mapsto \xi^{h_0}\sigma_{q;\xi}^{-h_1}W(\xi,m_1))(\tau^{h_2})=\sum_{n\ge h_0}W_{n-h_0}(m_1)q^{-h_1(n-h_0)}\frac{q^{\frac{n(n-1)}{2k}}}{q^{\frac{h_2n(h_2n-1)}{2k}}}\tau^{h_2 n}.$$
This yields
\begin{multline*}
\left|\hat{\mathcal{D}}_{h_2}(\xi\mapsto \xi^{h_0}\sigma_{q;\xi}^{-h_1}W(\xi,m_1))(\tau^{h_2})\right|\\
\le \sum_{n\ge h_0}\left[W_{n-h_0}(m_1)(1+|m_1|)^{\mu}e^{\beta|m_1|}\right](1+|m_1|)^{-\mu}e^{-\beta|m_1|}q^{-h_1(n-h_0)}\frac{q^{\frac{n(n-1)}{2k}}}{q^{\frac{h_2n(h_2n-1)}{2k}}}R^{h_2 n}\\
\le (1+|m_1|)^{-\mu}e^{-\beta|m_1|}\sum_{n\ge h_0}\left\|W_{n-h_0}(m)\right\|_{(\beta,\mu)} q^{-h_1(n-h_0)}\frac{q^{\frac{n(n-1)}{2k}}}{q^{\frac{h_2n(h_2n-1)}{2k}}}R^{h_2 n}\\
\le C_q(1+|m_1|)^{-\mu}e^{-\beta|m_1|}\sum_{n\ge h_0}\left\|W_{n-h_0}(m)\right\|_{(\beta,\mu)}R^{n}\\
\le C_q(1+|m_1|)^{-\mu}e^{-\beta|m_1|}R^{h_0}\left\|W(\tau,m)\right\|_{(1,R)}
\end{multline*}
Let $C_{\alpha,q}>0$ such that
$$\sup_{0<|\tau|<R}\frac{1}{|\tau|}\exp\left(-\frac{k}{2\log(q)}\log^2|\tau|-\alpha\log|\tau|\right)\le C_{\alpha,q},$$
The expression in (\ref{e421bb}) can be upper bounded by
\begin{multline*}C_{\alpha,q}C_3\left(\sum_{p\ge0}\left(\frac{R}{R_1}\right)^p\right)C_qR^{h_0}\left\|W(\tau,m)\right\|_{(1,R)}\left\|A(m)\right\|_{(\beta,\mu)}C_B\\
\times (1+|m|)^{\mu-\tilde{\alpha}}\int_{-\infty}^{\infty}\frac{1}{(1+|m-m_1|)^{\mu}(1+|m_1|)^{\mu-\hbox{deg}(B)}}dm_1\\
\le C_{A,B}C_{\alpha,q}C_3\left(\sum_{p\ge0}\left(\frac{R}{R_1}\right)^p\right)C_qR^{h_0}\left\|A(m)\right\|_{(\beta,\mu)}C_B\left\|W(\tau,m)\right\|_{(1,R)}.
\end{multline*}

Step 2.2: Assume that $\tau\in S_d$ with $|\tau|\ge R$. Observe that for $\xi\in\C$ with $|\xi|\le R$ one has that
$$\left|\xi^{h_0}\sigma_{q;\xi}^{-h_1}W(\xi,m)\right|\le R^{h_0}(1+|m|)^{-\mu}e^{-\beta|m|}\left\|W(\tau,m)\right\|_{(1,R)}.$$
In addition to this, taking into account Proposition~\ref{prop233} one has there exists $\Delta>0$ with
\begin{multline}\label{e422}
\left|L(\tau,m)\int_{-\infty}^{\infty}A(m-m_1)\hat{\mathcal{D}}_{h_2}\left(\xi\mapsto \xi^{h_0}\sigma_{q;\xi}^{-h_1}W(\xi,m_1)\right)(\tau^{h_2})B(im_1)dm_1\right|\\
\times (1+|m|)^{\mu}e^{\beta|m|}\frac{1}{|\tau|}\exp\left(-\frac{k}{2\log(q)}\log^2|\tau|-\alpha\log|\tau|\right)\\
\le C_1\exp\left(-\frac{\tilde{K}\log^2|\tau|}{2\log(q)}\right)|\tau|^h\frac{1}{(1+|m|)^{\tilde{\alpha}}}\hfill\\
\hfill\times K\exp\left(\frac{k}{h_2^2-1}\frac{1}{2\log(q)}\log^2(|\tau|^{h_2}+\Delta)+\alpha\log(|\tau|^{h_2}+\Delta)\right)\left\|A(m)\right\|_{(\beta,\mu)}R^{h_0}\\
\hfill\times\left\|W(\tau,m)\right\|_{(1,R)}\int_{-\infty}^{\infty}e^{-\beta|m-m_1|}\frac{1}{(1+|m-m_1|)^{\mu}}C_B(1+|m_1|)^{\hbox{deg}(B)-\mu}e^{-\beta|m_1|}dm_1\\
\times (1+|m|)^{\mu}e^{\beta|m|}\frac{1}{|\tau|}\exp\left(-\frac{k}{2\log(q)}\log^2|\tau|-\alpha\log|\tau|\right).
\end{multline}

Analogous estimates as above guarantee that the previous expression is upper bounded by
$$\mathcal{A}_2(|\tau|)C_1C_{A,B}KC_{B} \left\|A(m)\right\|_{(\beta,\mu)}R^{h_0}\left\|W(\tau,m)\right\|_{(1,R)},$$
with
\begin{multline*}
\mathcal{A}_2(|\tau|)=|\tau|^{h-1} \exp\left(-\frac{\tilde{K}\log^2|\tau|}{2\log(q)}\right)\exp\left(\frac{k}{h_2^2-1}\frac{1}{2\log(q)}\log^2(|\tau|^{h_2}+\Delta)+\alpha\log(|\tau|^{h_2}+\Delta)\right)\\
\times\exp\left(-\frac{k}{2\log(q)}\log^2|\tau|-\alpha\log|\tau|\right).
\end{multline*}
We observe that
$$\sup_{|\tau|\ge R}\mathcal{A}_2(|\tau|)\le \tilde{C}_{\mathcal{A}_2}\sup_{|\tau|\ge R}|\tau|^{N}\exp\left(\frac{1}{2\log(q)}\left(\frac{k}{h_2^2-1}-\tilde{K}\right)\log^2|\tau|\right)\le C_{\mathcal{A}_2},$$
for some $N\in\Z$ and some $\tilde{C}_{\mathcal{A}_2},C_{\mathcal{A}_2}>0$. The hypothesis made guarantee that the previous expression is upper bounded, on order to conclude the result for the constant
$$\tilde{C}_4=\tilde{C}_2+\left\|A(m)\right\|_{(\beta,\mu)}\max\left\{C_{\alpha,q}CC_BC_{A,B}C_qR^{h_0}\left(\sum_{p\ge 1}\left(\frac{R}{R_1}\right)^{p}\right),C_{\mathcal{A}_2}C_1C_{A,B}KC_{B} R^{h_0}\right\}.$$
\end{proof}

\section{Main problem}\label{secmainpr}

In this section, we state the main problem under study in the present work.

Let $q>1$ and let $k>0$. We also consider $\beta>0$ and $\mu>1$. We fix a finite set of tuples of nonnegative integers $\mathcal{A}\subseteq\N^3$ which satisfies that
\begin{equation}\label{ee1}
\ell_1\le \ell_0/k-1,\quad \underline{\ell}=(\ell_0,\ell_1,\ell_2)\in\mathcal{A}.
\end{equation} 
We also fix $\alpha_D>0$ and an integer number $d_{D}\ge1$ such that
\begin{equation}\label{e676c}
d_{D}>\max_{(\ell_0,\ell_1,\ell_2)\in\mathcal{A}}\left(\frac{k}{\ell_2^2-1}\right)^{1/2}.
\end{equation}

Let $Q(X),R_{D}(X),R_{\underline{\ell}}(X)\in\C[X]$ for all $\underline{\ell}\in\mathcal{A}$. We assume that 
$$Q(im)\neq 0,\qquad R_{D}(im)\neq 0,\qquad m\in\R.$$
Moreover, we assume the previous polynomials satisfy there exist $0<r_{Q,R_{D},1}<r_{Q,R_{D},2}$ with 
$$r_{Q,R_{D},1}\le\left|\frac{Q(im)}{R_{D}(im)}\right|\le r_{Q,R_{D},2},\qquad m\in\R.$$
In particular, one should have that $\hbox{deg}(Q)=\hbox{deg}(R_D)$. The constants $r_{Q,R_{D},1}$ and $r_{Q,R_{D},2}$ are small enough such that further conditions described in Section~\ref{secgeo} are fulfilled. We also assume that
\begin{equation}\label{e1910}
\hbox{deg}(R_{\underline{\ell}})\le \hbox{deg}(R_{D})=\hbox{deg}(Q),\qquad \underline{\ell}\in\mathcal{A}.
\end{equation}
It holds that
\begin{multline}\label{e461}
C_{R_D}(1+|m|)^{\hbox{deg}(R_D)}\le|R_{D}(im)|\le \tilde{C}_{R_D}(1+|m|)^{\hbox{deg}(R_D)},\\
 |R_{\underline{\ell}}(im)|\le C_{R_{\underline{\ell}}}(1+|m|)^{\hbox{deg}(R_{\underline{\ell}})},\quad |Q(im)|\ge C_{Q}(1+|m|)^{\hbox{deg}(Q)},
\end{multline}
for some constants $C_{R_D},\tilde{C}_{R_D},C_{R_{\underline{\ell}}},C_Q>0$, valid for all $m\in\R$.

We consider the problem
\begin{multline}
Q(\partial_z)u(t,z)=\exp_q\left(\alpha_Dt^{d_{D}}\sigma_{q;t}^{\frac{d_D}{k}}\right)R_D(\partial_z)u(t,z)+\sum_{\underline{\ell}=(\ell_0,\ell_1,\ell_2)\in\mathcal{A}}a_{\underline{\ell}}(z)\left(t^{\ell_0}\sigma_{q;t}^{\ell_1}R_{\underline{\ell}}(\partial_z)u\right)(t^{\ell_2},z)\\
+f(t,z)\label{epral}
\end{multline}
under vanishing initial data $u(0,z)\equiv 0$.

The function $x\mapsto \exp_q(x)$ consists on the $q$-analog of the exponential function whose definition and main properties are recalled in Section~\ref{secapqex1}. We define the formal operator 
$$\exp_q(\alpha_Dt^{d_D}\sigma_{q;t}^{\frac{d_D}{k}})=\sum_{n\ge0}\frac{1}{[n]_q^{!}}\alpha_D^{n}(t^{d_D}\sigma_{q;t}^{\frac{d_D}{k}})^{(n)},$$
where $(t^{d_D}\sigma_{q;t}^{\frac{d_D}{k}})^{(n)}$ stands for the $n$-th iteration of the $q$-difference operator $t^{d_D}\sigma_{q;t}^{\frac{d_D}{k}}$.

The coefficients $a_{\underline{\ell}}(z)$ are constructed as follows. For all $\underline{\ell}\in\mathcal{A}$, we consider a continuous function $A_{\underline{\ell}}:\R\to\C$ such that there exists $C_{\ell}>0$ with
\begin{equation}\label{e588}
|A_{\underline{\ell}}(m)|\le C_{\underline{\ell}}\frac{1}{(1+|m|)^{\mu}}\exp(-\beta|m|),\quad m\in\R.
\end{equation}
We define $a_{\underline{\ell}}(z)=\mathcal{F}^{-1}(A_{\underline{\ell}})(z)$ for all $\underline{\ell}\in \mathcal{A}$. The function $a_{\underline{\ell}}\in\mathcal{O}(H_{\beta})$, with 
$$H_{\beta}=\{z\in\C:|\hbox{Im}(z)|<\beta\},$$
in view of Section~\ref{ap2}.

The forcing term $f$ is built as follows. Let $J\subseteq\N$ be a finite set, and choose $d_1\in\R$. For every $j\in J$ we consider a continuous function $\mathcal{F}_j:\R\to\C$ such that there exists $C_{j}>0$ with
\begin{equation}\label{e596}
|\mathcal{F}_{j}(m)|\le C_{j}\frac{1}{(1+|m|)^{\mu}}\exp(-\beta|m|),\quad m\in\R.
\end{equation}
Let $\mathcal{F}:\C\times\R\to\C$ be defined by
$$\mathcal{F}(\tau,m)=\sum_{j\in J}\mathcal{F}_j(m)\tau^j.$$
We fix 
$$f(t,z)=\frac{\pi_{q,k}}{(2\pi)^{1/2}}\int_{L_{d_1}}\int_{-\infty}^{\infty}\mathcal{F}(\tau,m)\Theta_k\left(\frac{t}{\tau}\right)e^{izm}\frac{d\tau}{\tau}dm,$$
with 
$$\pi_{q,k}:=\frac{q^{-\frac{1}{8k}}\sqrt{k}}{(2\pi \log(q))^{1/2}},$$
and where the integration path is $L_{d_1}=[0,\infty)e^{id_1}$. The kernel function $\Theta_k$ is defined by 
$$\Theta_k(z):=q^{-\frac{k}{2}\frac{\log(z)}{\log(q)}\left(\frac{\log(z)}{\log(q)}-\frac{1}{k}\right)}.$$

Observe from Section~\ref{sec200} and Section~\ref{ap2} that $f$ is the $q$-Laplace and inverse Fourier transform of $\mathcal{F}(\tau,m)$. Therefore, $f(t,z)$ turns out to be analytic on $\tilde{D}(0,R)\times H_{\beta'}$, where $\tilde{D}(0,R)$ stands for the universal covering of $D(0,R)$, for some $R>0$ and $0<\beta'<\beta$. Moreover, in view of Proposition~\ref{prop254}, one has that
$$f(t,z)=\sum_{j\in J}F_j(z)q^{\frac{j(j-1)}{2k}}t^j,\quad (t,z)\in \tilde{D}(0,R)\times H_{\beta'}, $$
where
$$F_j(z)=\frac{1}{(2\pi)^{1/2}}\int_{-\infty}^{\infty}\mathcal{F}_j(m)e^{izm}dm,\quad z\in H_{\beta'}.$$

\subsection{Solution strategy}\label{secI}

The strategy to solve (\ref{epral}) is to search for solutions in the form of an inverse Fourier transform, i.e.
$$u(t,z)=\frac{1}{(2\pi)^{1/2}}\int_{-\infty}^{\infty}U(t,m)e^{izm}dm,$$
for some expression $U(t,m)$. In view of the properties of inverse Fourier transform in Section~\ref{ap2}, we derive that $u(t,z)$ formally solves (\ref{epral}) if $U(t,m)$ is a solution of the following integro-$q$-difference Mahler equation.
\begin{multline}
Q(im)U(t,m)=\exp_q\left(\alpha_D t^{d_D}\sigma_{q;t}^{\frac{d_D}{k}}\right)R_D(im)U(t,m)\\
+\sum_{\underline{\ell}=(\ell_0,\ell_1,\ell_2)\in\mathcal{A}}\frac{1}{(2\pi)^{1/2}}\int_{-\infty}^{+\infty}A_{\underline{\ell}}(m-m_1)(t^{\ell_0}\sigma_{q;t}^{\ell_1}U)(t^{\ell_2},m_1)R_{\underline{\ell}}(im_1)dm_1+\sum_{j\in J}\mathcal{F}_j(m)q^{\frac{j(j-1)}{2k}}t^j,\label{e169}
\end{multline}
for vanishing initial data $U(0,m)\equiv 0$.

In the next step, we seek for solutions to (\ref{e169}) in the form of a formal power series
\begin{equation}\label{e173}
\hat{U}(t,m)=\sum_{n\ge1}U_{n}(m)t^{n}.
\end{equation}
We ask the coefficients $m\mapsto U_n(m)$ to satisfy bounds in the form
$$|U_n(m)|\le C\frac{1}{(1+|m|)^{\mu}}\exp(-\beta|m|),\quad m\in\R,$$
for some $C>0$, common for all $n$. In the terminology of Section~\ref{ap2}, the functions $m\mapsto U_n(m)$ belong to $E_{(\beta,\mu)}$ for $\beta>0$ and $\mu>1$, with a common upper bound for their norm in such Banach space.

In view of the results of Section~\ref{secb}, we proceed with a second step in our strategy to solve the main problem under study, after the first step, leading to the auxiliary equation (\ref{e169}).

Let us denote the formal $q$-Borel transform of order $k$ of the formal series $t\mapsto \hat{U}(t,m)$, defined in (\ref{e713}), by
$$W(\tau,m)=\hat{\mathcal{B}}_{q;k}(t\mapsto \hat{U}(t,m))(\tau),$$
for all $m\in\R$.

According to Proposition~\ref{prop195}, we have that
\begin{equation}\label{e286a}
\hat{\mathcal{B}}_{q;k}\left(\exp_q\left(\alpha_D t^{d_D}\sigma_{q;t}^{\frac{d_D}{k}}\right)R_D(im)\hat{U}(t,m)\right)(\tau)=\exp_q\left(\alpha_D\frac{\tau^{d_{D}}}{q^{\frac{d_D(d_{D}-1)}{2k}}}\right)R_D(im)W(\tau,m).
\end{equation}
In addition to this, if $\ell=(\ell_0,\ell_1,\ell_2)\in\mathcal{A}$ with $\ell_2=1$, then Proposition~\ref{prop195} leads to
\begin{equation}\label{e286b}
\hat{\mathcal{B}}_{q;k}\left( t^{\ell_0}\sigma_{q;t}^{\ell_1}\hat{U}(t,m_1)\right)(\tau)=\frac{\tau^{\ell_0}}{q^{\frac{\ell_0(\ell_0-1)}{2k}}}\sigma_{q;\tau}^{\ell_1-\frac{\ell_0}{k}}W(\tau,m_1).
\end{equation}
On the other hand, if $\ell=(\ell_0,\ell_1,\ell_2)\in\mathcal{A}$ with $\ell_2\ge2$, in view of Proposition~\ref{prop232} and Proposition~\ref{prop233} we get that
\begin{multline}
\hat{\mathcal{B}}_{q;k}\left( (t^{\ell_0}\sigma_{q;t}^{\ell_1}\hat{u})(t^{\ell_2},m_1)\right)(\tau)\\
=\Xi_{q,k'}\int_{\tilde{C}_{\rho/2}}\mathbb{D}_{q,k'}(x/\tau^{\ell_2})\sigma_{q;x}^{-k''}\left[\frac{x^{\ell_0}}{q^{\frac{\ell_0(\ell_0-1)}{2k}}}\sigma_{q;x}^{\ell_1-\frac{\ell_0}{k}}W(x,m_1)\right]\frac{dx}{x},\label{e286c}
\end{multline}
where $\tilde{C}_{\rho/2}=\left\{\frac{\rho}{2}e^{i\theta}:\theta\in\R\right\}\subseteq\tilde{D}(0,R)$, $\Xi_{q,k'}:=-i\frac{q^{1/(8k')}\sqrt{k'}}{(2\pi \log(q))^{1/2}}$, with $k'=\frac{k}{\ell_2^2-1}$ and $k''=\frac{\ell_2^2-\ell_2}{2k}$, together with
$$\mathbb{D}_{q,k'}(x/\tau^{\ell_2})=q^{\frac{k'}{2}\frac{\log(x/\tau^{\ell_2})}{\log(q)}\left(\frac{\log(x/\tau^{\ell_2})}{\log(q)}-\frac{1}{k'}\right)}.$$

In view of (\ref{e286a}), (\ref{e286b}) and (\ref{e286c}), we get that the formal power series $\hat{U}(t,m)$ in (\ref{e173}) is a solution to (\ref{e169}) if the formal $q$-Borel transform 
$$W(\tau,m)=\hat{\mathcal{B}}_{q;k}(\hat{U})(\tau)$$
satisfies the next $q$-difference-integral equation.

\begin{multline}\label{eaux1}
W(\tau,m)=\sum_{\underline{\ell}=(\ell_0,\ell_1,\ell_2)\in\mathcal{A},\ell_2=1}\frac{1}{(2\pi)^{1/2}}\frac{1}{P_m(\tau)}\int_{-\infty}^{\infty}A_{\underline{\ell}}(m-m_1)\frac{\tau^{\ell_0}}{q^{\frac{\ell_0(\ell_0-1)}{2k}}}\sigma_{q;\tau}^{\ell_1-\frac{\ell_0}{k}}W(\tau,m_1)R_{\underline{\ell}}(im_1)dm_1\\
+\sum_{\underline{\ell}=(\ell_0,\ell_1,\ell_2)\in\mathcal{A},\ell_2\ge 2}\frac{1}{(2\pi)^{1/2}}\frac{1}{P_m(\tau)}\int_{-\infty}^{\infty}A_{\underline{\ell}}(m-m_1)\hat{\mathcal{D}}_{\ell_2}\left(\xi\mapsto\frac{\xi^{\ell_0}}{q^{\frac{\ell_0(\ell_0-1)}{2k}}}\sigma_{q;\xi}^{\ell_1-\frac{\ell_0}{k}}W(\xi,m_1)\right)(\tau^{\ell_2})\\
\times R_{\underline{\ell}}(im_1)dm_1+\sum_{j\in J}\mathcal{F}_j(m)\frac{\tau^j}{P_m(\tau)},
\end{multline}
where
$$P_m(\tau)=Q(im)-\exp_q(\tilde{\alpha}_D\tau^{d_D})R_D(im),$$
with 
\begin{equation}\label{e382}
\tilde{\alpha}_D=\frac{\alpha_D}{q^{\frac{d_D(d_D-1)}{2k}}}.
\end{equation}

In view of the second statement of Proposition~\ref{prop233}, one may also consider the functional equation

\begin{multline}\label{eaux2}
W(\tau,m)=\sum_{\underline{\ell}=(\ell_0,\ell_1,\ell_2)\in\mathcal{A},\ell_2=1}\frac{1}{(2\pi)^{1/2}}\frac{1}{P_m(\tau)}\int_{-\infty}^{\infty}A_{\underline{\ell}}(m-m_1)\frac{\tau^{\ell_0}}{q^{\frac{\ell_0(\ell_0-1)}{2k}}}\sigma_{q;\tau}^{\ell_1-\frac{\ell_0}{k}}W(\tau,m_1)R_{\underline{\ell}}(im_1)dm_1\\
+\sum_{\underline{\ell}=(\ell_0,\ell_1,\ell_2)\in\mathcal{A},\ell_2\ge2}\frac{1}{(2\pi)^{1/2}}\frac{1}{P_m(\tau)}\int_{-\infty}^{\infty}A_{\underline{\ell}}(m-m_1)\Xi_{q,k'}\int_{\tilde{C}_{\rho/2}}\mathbb{D}_{q,k'}(x/\tau^{\ell_2})\\
\times\sigma_{q;x}^{-k''}\left[\frac{x^{\ell_0}}{q^{\frac{\ell_0(\ell_0-1)}{2k}}}\sigma_{q;x}^{\ell_1-\frac{\ell_0}{k}}W(x,m_1)\right]\frac{dx}{x} R_{\underline{\ell}}(im_1)dm_1+\sum_{j\in J}\mathcal{F}_j(m)\frac{\tau^j}{P_m(\tau)},
\end{multline}

with $k'=\frac{k}{\ell_2^2-1}$ and $k''=\frac{\ell_2^2-\ell_2}{2k}$.

\subsection{Geometric configuration for the solution and auxiliary estimates}\label{secgeo}

In this section, we fix the geometric configuration of some elements involved in the construction of the main problem under study. More precisely, we select an unbounded sector of bisecting direction $d\in\R$, say $S_d$, and choose $\rho>0$ such that the following statements hold, regarding the unbounded sector $U$ described in 	Proposition~\ref{prop430}:
\begin{itemize}
\item $\tilde{\alpha}_{D}\tau^{d_{D}}\in U$ for all $\tau\in S_d$, where $\tilde{\alpha}_D$ is given in (\ref{e382}),
\item $\tilde{\alpha}_{D}\tau^{d_{D}}\in D(0,\frac{q^{1/2}}{q-1})$ for all $\tau\in D(0,\rho)$.
\end{itemize}

Observe that the previous constraints are feasible under an appropriate choice of $d\in\R$ and with $\rho$ given by $0<\rho\le\left(\frac{q^{1/2}}{q-1}\frac{q^{\frac{d_D(d_D-1)}{2k}}}{\alpha_D}\right)^{1/d_{D}}$.

Let us denote for every $x>0$
$$\mu(x):=\frac{\log^2(x)}{2\log(q)}+\left(-\frac{1}{2}+\frac{\log(q-1)}{\log(q)}\right)\log(x).$$

As a result of Proposition~\ref{prop430}, we deduce that for all $\tau\in S_d$ with $|\tau|\ge \rho$ then one has that
\begin{equation}\label{e398}
|\exp_q(\tilde{\alpha}_D\tau^{d_{D}})|\le K_1\exp(\mu(\tilde{\alpha}_D|\tau|^{d_{D}})),
\end{equation}
together with
\begin{equation}\label{e398b}
|\exp_q(\tilde{\alpha}_D\tau^{d_{D}})|\ge \frac{\epsilon}{K_0}\exp(\mu(\tilde{\alpha}_D|\tau|^{d_{D}})),
\end{equation}
for some well chosen constants $\epsilon,K_0,K_1>0$. Observe that 
$$\exp(\mu(\tilde{\alpha}_D|\tau|^{d_{D}}))\ge C_{\mu}\exp\left(\frac{d_D^2\log^2|\tau|}{2\log(q)}\right)|\tau|^{H},$$
for some $C_{\mu}>0$ and some $H\in\R$, for all $\tau\in S_d$ with $|\tau|\ge \rho$.

For $\tau\in \overline{D}(0,\rho)$ one also has that
\begin{equation}\label{e398c}
|\exp_q(\tilde{\alpha}_D\tau^{d_{D}})|\ge C_0,
\end{equation}
for a suitable constant $C_0>0$.

We choose $r_{Q,R_{D},2}$ in the statement of the problem (see Section~\ref{secmainpr}) such that 
\begin{equation}\label{e412}
r_{Q,R_{D},2}<\min\{C_0,\frac{\epsilon}{K_0}\min_{x\ge \frac{q^{1/2}}{q-1}}e^{\mu(x)}\}.
\end{equation}

We now provide estimates from above for some expression associated to the geometric configuration of the problem. Let us define for all $m\in\R$ and $\tau\in\C$ 
\begin{equation}\label{e419}
P_m(\tau):=Q(im)-\exp_q(\tilde{\alpha}_D\tau^{d_D})R_D(im).
\end{equation}
By construction, we observe there exists $\delta_1>0$ such that
$$\left|\frac{Q(im)}{R_{D}(im)}-\exp_q(\tilde{\alpha}_D\tau^{d_D})\right|\ge \delta_1,\quad \tau\in S_d\cup D(0,\rho),\quad m\in\R.$$
As a result we get the lower bounds
\begin{equation}\label{e421}
|P_m(\tau)|\ge \delta_1|R_D(im)|,\quad \tau\in S_d\cup D(0,\rho),\quad m\in\R.
\end{equation}

We need sharper estimates for our purposes when considering large values of $|\tau|$. Indeed, one can factorize 
$$\exp_q(\tilde{\alpha}_D\tau^{d_D})-\frac{Q(im)}{R_{D}(im)}=\exp_q(\tilde{\alpha}_D\tau^{d_D})\left(1-\frac{Q(im)}{R_{D}(im)}\frac{1}{\exp_q(\tilde{\alpha}_D\tau^{d_D})}\right).$$
From the choice of $U$ at the beginning of this section, we derive that
$$|\exp_q(\tilde{\alpha}_D\tau^{d_D})|\ge \frac{\epsilon}{K_0}\exp(\mu(\tilde{\alpha}_D|\tau|^{d_D})),$$
for every $\tau\in S_d$ with $|\tau|\ge \rho$. In addition to this, and due to the fact that $\mu(x)\to+\infty$ when $x\to+\infty$, we arrive at the existence of $r_1\ge 0$ such that if $\tau\in S_d$ with $|\tau|\ge r_1$ and all $m\in \R$, then one has that
$$\left|\frac{Q(im)}{R_{D}(im)}\right|\frac{1}{|\exp_q(\tilde{\alpha}_D\tau^{d_D})|}\le r_{Q,R_D,2}\frac{K_0}{\epsilon}\frac{1}{\exp(\mu(\tilde{\alpha}_D|\tau|^{d_D}))}\le\frac{1}{2}.$$
This entails that 
$$\left|\exp_q(\tilde{\alpha}_D\tau^{d_D})\left(1-\frac{Q(im)}{R_{D}(im)}\frac{1}{\exp_q(\tilde{\alpha}_D\tau^{d_D})}\right)\right|\ge|\exp_q(\tilde{\alpha}_D\tau^{d_D})|\frac{1}{2},$$
for all $\tau\in S_d$ with $|\tau|\ge r_1$ and all $m\in \R$. A compactness argument allows to modify the value of $r_1$ adequately, providing analogous bounds. 

As a result, we have achieved the following result.

\begin{lemma}\label{lema437}
There exists $\delta_1>0$ such that
$$|P_m(\tau)|\ge \delta_1|R_D(im)|,\quad \tau\in S_d\cup D(0,\rho),\quad m\in\R.$$
Moreover, for every $r_1>0$ there exists $M>0$ such that
$$|P_m(\tau)|\ge \frac{\epsilon}{MK_0}\exp(\mu(\tilde{\alpha}_D|\tau|^{d_D})) |R_D(im)|,$$
for every $\tau\in S_d$ with $|\tau|\ge r_1$ and all $m\in\R$, or in other words,
$$|P_m(\tau)|\ge \frac{\epsilon C_{\mu}}{MK_0}\exp(\frac{d_D^2\log^2|\tau|}{2\log(q)})|\tau|^H |R_D(im)|,$$
for $C_{\mu},H$ as described after the definition of $\mu(\cdot)$.
\end{lemma}




The next result will be crucial when determining the solution of the auxiliary problem on some neighborhood of the origin.

\begin{lemma}\label{lema644}
Let us write $1/P_m(\tau)=\sum_{p\ge0}f_p(m)\tau^p$, for all $m\in\R$ and $\tau\in D(0,\rho)$. Then, for every $0<R_1<\rho$ there exists $C_{P}>0$ such that 
$$|f_p(m)|\le C_P\frac{1}{R_1^{p}}\frac{1}{(1+|m|)^{\hbox{deg}(R_D)}},\qquad m\in\R.$$
\end{lemma}
\begin{proof}
Let $m\in\R$, $p\ge0$ be an integer, and $0<R_1<\rho$. A direct application of Cauchy's integral formula yields
$$|f_p(m)|=\left|\partial_{\tau}^{p}\left(\frac{1}{P_m(\tau)}\right)_{\tau=0}\right|\frac{1}{p!}\le\sup_{|\omega|=R_1}\left|\frac{1}{P_m(\omega)}\right|\frac{1}{R_1^p}\\
\le \frac{1}{\delta_1C_{R_D}(1+|m|)^{\hbox{deg}(R_D)}}\frac{1}{R_1^p}.$$
From Lemma~\ref{lema437}, we get a constant $\delta_1>0$ such that
$$\sup_{|\omega|=R_1}\frac{1}{|P_{m}(\omega)|}\le \frac{1}{\delta_1}\frac{1}{|R_{D}(im)|}$$
for all $m\in\R$. Besides, we have applied the fact that $R_{D}(im)\neq 0$ for all $m\in\R$, which guarantees the existence of $C_{R_D}>0$ with $|R_{D}(im)|\ge C_{R_D}(1+|m|)^{\hbox{deg}(R_D)}$.
\end{proof}

\section{Preliminary results}

In this section, we determine preliminary auxiliary results to be considered in the construction of the analytic solution of the main problem. This is performed by a fixed point argument involving the Banach space of functions detailed in Section~\ref{sec22} and Section~\ref{sec23}, respectively. We maintain all the constructions and assumptions made on the main elements involving the main problem under study, equation (\ref{epral}), in Section~\ref{secmainpr}.

\subsection{Fixed point in a Banach space, I}\label{secbsi}

In this section, we guarantee the existence of a fixed point for certain operator defined in the Banach space developed in Section~\ref{sec22}. This fixed point gives rise to a solution of the auxiliary problem (\ref{eaux1}), under the geometric configuration described in Section~\ref{secgeo}.

Let us consider the operator 
\begin{multline}\label{e614}
\mathcal{H}_1(\omega(\tau,m))=\sum_{\underline{\ell}=(\ell_0,\ell_1,\ell_2)\in\mathcal{A},\ell_2=1}\frac{1}{(2\pi)^{1/2}}\frac{1}{P_m(\tau)}\int_{-\infty}^{\infty}A_{\underline{\ell}}(m-m_1)\frac{\tau^{\ell_0}}{q^{\frac{\ell_0(\ell_0-1)}{2k}}}\sigma_{q;\tau}^{\ell_1-\frac{\ell_0}{k}}\omega(\tau,m_1)R_{\underline{\ell}}(im_1)dm_1\\
+\sum_{\underline{\ell}=(\ell_0,\ell_1,\ell_2)\in\mathcal{A},\ell_2\ge 2}\frac{1}{(2\pi)^{1/2}}\frac{1}{P_m(\tau)}\int_{-\infty}^{\infty}A_{\underline{\ell}}(m-m_1)\hat{\mathcal{D}}_{\ell_2}\left(\xi\mapsto\frac{\xi^{\ell_0}}{q^{\frac{\ell_0(\ell_0-1)}{2k}}}\sigma_{q;\xi}^{\ell_1-\frac{\ell_0}{k}}\omega(\xi,m_1)\right)(\tau^{\ell_2})\\
\times R_{\underline{\ell}}(im_1)dm_1+\sum_{j\in J}\mathcal{F}_j(m)\frac{\tau^j}{P_m(\tau)}.
\end{multline}

\begin{prop}\label{prop808}
Let $\rho>0$ and $\varpi>0$. Then, there exists $0<R<\rho$ such that the operator $\mathcal{H}_1$ defined in (\ref{e614}) admits a unique fixed point $\omega_1(\tau,m)\in HE_{(\beta,\mu,R)}$, with $\left\|\omega_1(\tau,m)\right\|_{(1,R)}\le\varpi$.
\end{prop}
\begin{proof}

Let $\varpi>0$ and take $\omega(\tau,m)\in HE_{(\beta,\mu,R)}$ with $\left\|\omega(\tau,m)\right\|_{(1,R)}\le\varpi$. Under the condition (\ref{e1910}) and the assumption $\ell_0\ge \ell_1 k$, we apply Lemma~\ref{lema3} and Lemma~\ref{lema644} to arrive at a constant $\tilde{C}_1>0$ with
$$\left\|\frac{1}{P_m(\tau)}\int_{-\infty}^{\infty}A_{\underline{\ell}}(m-m_1)\tau^{\ell_0}\sigma_{q;\tau}^{\ell_1-\frac{\ell_0}{k}}\omega(\tau,m_1)R_{\underline{\ell}}(im_1)dm_1\right\|_{(1,R)}\le \tilde{C}_1\left\|\omega(\tau,m)\right\|_{(1,R)},$$
provided that the radius $R>0$ fulfills $0<R<R_1$ for $R_1$ appearing in Lemma~\ref{lema644}.

The application of Lemma~\ref{lema4} and Lemma~\ref{lema644} under the condition (\ref{e1910}) and the assumption (\ref{ee1}) yields a constant $\tilde{C}_2>0$ with

\begin{multline*}
\left\|\frac{1}{P_m(\tau)}\int_{-\infty}^{\infty}A_{\underline{\ell}}(m-m_1)\hat{\mathcal{D}}_{\ell_2}\left(\xi\mapsto\xi^{\ell_0}\sigma_{q;\xi}^{\ell_1-\frac{\ell_0}{k}}\omega(\xi,m_1)\right)(\tau^{\ell_2}) R_{\underline{\ell}}(im_1)dm_1\right\|_{(1,R)}\\
\le \tilde{C}_2\left\|\omega(\tau,m)\right\|_{(1,R)}.
\end{multline*}

Finally, in view of Lemma~\ref{lema644} and the construction of the forcing term in the main equation, we derive that for every $0<R<R_1<\rho$ and $j\in J$, one has
\begin{multline*}
\left\|\mathcal{F}_j(m)\frac{\tau^j}{P_m(\tau)}\right\|_{(1,R)}=\left\|\sum_{p\ge0}f_p(m)\mathcal{F}_j(m)\tau^{j+p}\right\|_{(1,R)}=\sum_{p\ge0}\left\|f_p(m)\mathcal{F}_j(m)\right\|_{(\beta,\mu)}R^{j+p}\\
\le\sum_{p\ge0}\sup_{m\in\R}|\mathcal{F}_j(m)f_p(m)|(1+|m|)^{\mu}e^{\beta|m|}R^{j+p}\le C_j\sum_{p\ge0}\sup_{m\in\R}|f_p(m)|R^{j+p}\le C_jC_PR^j\sum_{p\ge0}\left(\frac{R}{R_1}\right)^{p}.
\end{multline*}

Let $R>0$ be small enough such that
\begin{multline*}
\sum_{\underline{\ell}=(\ell_0,\ell_1,\ell_2)\in\mathcal{A},\ell_2=1}\tilde{C}_1\frac{1}{(2\pi)^{1/2}}\frac{1}{q^{\frac{\ell_0(\ell_0-1)}{2k}}}\\
+\sum_{\underline{\ell}=(\ell_0,\ell_1,\ell_2)\in\mathcal{A},\ell_2\ge 2}\tilde{C}_2\frac{1}{(2\pi)^{1/2}}\frac{1}{q^{\frac{\ell_0(\ell_0-1)}{2k}}}+C_P\sum_{j\in J}C_jR^j\sum_{p\ge0}\left(\frac{R}{R_1}\right)^{p}\le 1.
\end{multline*}

Observe that $\tilde{C}_1$ and $\tilde{C}_2$ can be chosen small enough with small enough $R>0$. See Remark at page~\pageref{remark}.

Then, one has that $\mathcal{H}_1(\omega(\tau,m))\in HE_{(\beta,\mu,R)}$, with $\left\|\mathcal{H}_1(\omega(\tau,m))\right|_{(1,R)}\le\varpi$.

We now consider $\omega_1(\tau,m),\omega_2(\tau,m)\in HE_{(\beta,\mu,R)}$ with $\left\|\omega_j(\tau,m)\right\|_{(1,R)}\le \varpi$ for $j=1,2$. Analogous estimates as above lead us to 

\begin{multline*}
\left\|\frac{1}{P_m(\tau)}\int_{-\infty}^{\infty}A_{\underline{\ell}}(m-m_1)\tau^{\ell_0}\sigma_{q;\tau}^{\ell_1-\frac{\ell_0}{k}}(\omega_1(\tau,m_1)-\omega_2(\tau,m_1))R_{\underline{\ell}}(im_1)dm_1\right\|_{(1,R)}\\
\le \tilde{C}_1\left\|\omega_1(\tau,m)-\omega_2(\tau,m)\right\|_{(1,R)},
\end{multline*}
together with
\begin{multline*}
\left\|\frac{1}{P_m(\tau)}\int_{-\infty}^{\infty}A_{\underline{\ell}}(m-m_1)\hat{\mathcal{D}}_{\ell_2}\left(\xi\mapsto\xi^{\ell_0}\sigma_{q;\xi}^{\ell_1-\frac{\ell_0}{k}}(\omega_1(\xi,m_1)-\omega_2(\xi,m_1))\right)(\tau^{\ell_2}) R_{\underline{\ell}}(im_1)dm_1\right\|_{(1,R)}\\
\le \tilde{C}_2\left\|\omega_1(\tau,m)-\omega_2(\tau,m)\right\|_{(1,R)}.
\end{multline*}

Given $R>0$ such that
$$\sum_{\underline{\ell}=(\ell_0,\ell_1,\ell_2)\in\mathcal{A},\ell_2=1}\tilde{C}_1\frac{1}{(2\pi)^{1/2}}\frac{1}{q^{\frac{\ell_0(\ell_0-1)}{2k}}}+\sum_{\underline{\ell}=(\ell_0,\ell_1,\ell_2)\in\mathcal{A},\ell_2\ge 2}\tilde{C}_2\frac{1}{(2\pi)^{1/2}}\frac{1}{q^{\frac{\ell_0(\ell_0-1)}{2k}}}\le \frac{1}{2},$$
we conclude that 
$$\left\|\mathcal{H}_1(\omega_1(\tau,m))-\mathcal{H}_1(\omega_2(\tau,m))\right\|_{(1,R)}\le \frac{1}{2}\left\|\omega_1(\tau,m)-\omega_2(\tau,m)\right\|_{(1,R)}.$$

The classical contractive mapping theorem applied to 
$$\mathcal{H}_1:\overline{D}(0,\varpi)\to \overline{D}(0,\varpi),$$
with $\overline{D}(0,\varpi)\subseteq HE_{(\beta,\mu,R)}$ allows us to conclude the existence of a unique fixed point of $\mathcal{H}_1$.
\end{proof}

\subsection{Fixed point in a Banach space, II}\label{secbsii}

In this subsection, we study the action of a linear integral operator on the Banach space of functions introduced in Definition~\ref{defi111}. 

In the whole section, we fix an unbounded sector with vertex at the origin and bisecting direction $d\in\R$, say $S_d$ as described in Section~\ref{secgeo}. We also set $r>0$.

Let us study the operator 

\begin{multline}
\mathcal{H}_2(\omega(\tau,m))=\sum_{\underline{\ell}=(\ell_0,\ell_1,\ell_2)\in\mathcal{A},\ell_2=1}\frac{1}{(2\pi)^{1/2}}\frac{1}{P_m(\tau)}\int_{-\infty}^{\infty}A_{\underline{\ell}}(m-m_1)\frac{\tau^{\ell_0}}{q^{\frac{\ell_0(\ell_0-1)}{2k}}}\sigma_{q;\tau}^{\ell_1-\frac{\ell_0}{k}}\omega(\tau,m_1)R_{\underline{\ell}}(im_1)dm_1\\
+\sum_{\underline{\ell}=(\ell_0,\ell_1,\ell_2)\in\mathcal{A},\ell_2\ge2}\frac{1}{(2\pi)^{1/2}}\frac{1}{P_m(\tau)}\int_{-\infty}^{\infty}A_{\underline{\ell}}(m-m_1)\Xi_{q,k'}\int_{\tilde{C}_{\rho/2}}\mathbb{D}_{q,k'}(x/\tau^{\ell_2})\\
\times\sigma_{q;x}^{-k''}\left[\frac{x^{\ell_0}}{q^{\frac{\ell_0(\ell_0-1)}{2k}}}\sigma_{q;x}^{\ell_1-\frac{\ell_0}{k}}\omega(x,m_1)\right]\frac{dx}{x} R_{\underline{\ell}}(im_1)dm_1+\sum_{j\in J}\mathcal{F}_j(m)\frac{\tau^j}{P_m(\tau)},\label{e715}
\end{multline}
with $k'=k/(\ell_2^2-1)$ and $k''=(\ell_2^2-\ell_2)/2k$.

\begin{prop}\label{prop882}
Let $\varpi>0$ and take $\omega(\tau,m)\in E^{d}_{(R,k,\alpha,\beta,\mu)}$ with $\left\|\omega(\tau,m)\right\|_{(R,k,\alpha,\beta,\mu)}\le \varpi$. Under the condition (\ref{e1910}) and the assumption (\ref{ee1}) for every $(\ell_0,\ell_1,\ell_2)\in\mathcal{A}$ with $\ell_2=1$, we apply Lemma~\ref{lema437} and Lemma~\ref{lema644} together with Lemma~\ref{lema486a}; one arrives at the existence of $\varsigma_{c},\varsigma_{f}>0$ such that if
$$C_{\underline{\ell}}\le \varsigma_{c},\qquad C_j\le \varsigma_{f},$$
for all $\underline{\ell}\in\mathcal{A}$ and all $j\in J$, and provided that $R>0$ is taken close enough to 0, the operator $\mathcal{H}_2$ defined in (\ref{e715}) admits a unique fixed point $\omega_2(\tau,m)\in E^d_{(R,k,\alpha,\beta,\mu)}$, with $\left\|\omega_2(\tau,m)\right\|_{(R,k,\alpha,\beta,\mu)}\le \varpi$.
\end{prop}

\begin{proof}
Let $\varpi>0$ and take $\omega(\tau,m)\in E^{d}_{(R,k,\alpha,\beta,\mu)}$ with 
$$\left\|\omega(\tau,m)\right\|_{(R,k,\alpha,\beta,\mu)}\le \varpi.$$
For every $(\ell_0,\ell_1,\ell_2)\in\mathcal{A}$ with $\ell_2=1$ we apply Lemma~\ref{lema437}, together with Lemma~\ref{lema486a}, one arrives at

$$\left\|\frac{1}{P_m(\tau)}\int_{-\infty}^{\infty}A_{\underline{\ell}}(m-m_1)\tau^{\ell_0}\sigma_{q;\tau}^{\ell_1-\frac{\ell_0}{k}}\omega(\tau,m_1)R_{\underline{\ell}}(im_1)dm_1\right\|_{(R,k,\alpha,\beta,\mu)}\le \tilde{C}_3\left\|\omega(\tau,m)\right\|_{(R,k,\alpha,\beta,\mu)},$$
for $\tilde{C}_3>0$ determined in Lemma~\ref{lema486a}. Observe from the proof of Lemma~\ref{lema486a} that $\tilde{C}_3$ can be arbitrary small by choosing $\tilde{C}_1$ and $\left\|A_{\underline{\ell}}(m)\right\|_{(\beta,\mu)}$ small enough. The constant $\tilde{C}_1>0$ can be diminished by choosing $R>0$ small enough, whereas $\left\|A_{\underline{\ell}}(m)\right\|_{(\beta,\mu)}$ coincides with $C_{\underline{\ell}}$ for all $\underline{\ell}\in\mathcal{A}$.

Provided that conditions (\ref{e1910}) and  (\ref{e676c}), and the assumption (\ref{ee1}) applied to $\underline{\ell}=(\ell_0,\ell_1,\ell_2)\in\mathcal{A}$ with $\ell_2\ge 2$ hold, taking into account Lemma~\ref{lema437} and Lemma~\ref{lema644}, and the representation in Proposition~\ref{prop233}, one arrives at
\begin{multline*}
\left\|\frac{1}{P_m(\tau)}\int_{-\infty}^{\infty}A_{\underline{\ell}}(m-m_1)\int_{\tilde{C}_{\rho/2}}\mathbb{D}_{q,k'}(x/\tau^{\ell_2})\sigma_{q;x}^{-k''}\left[x^{\ell_0}\sigma_{q;x}^{\ell_1-\frac{\ell_0}{k}}\omega(x,m_1)\right]\frac{dx}{x} R_{\underline{\ell}}(im_1)dm_1\right\|_{(R,k,\alpha,\beta,\mu)}\\
\le \tilde{C}_4 \left\|\omega(\tau,m)\right\|_{(R,k,\alpha,\beta,\mu)},
\end{multline*}
for $\tilde{C}_4$ determined in Lemma~\ref{lema486b}. In particular, observe that $\tilde{C}_4$ tends to zero with $\tilde{C}_2$ (and therefore $R$, regarding the proof of Lemma~\ref{lema4}) and by choosing small enough $R$.

Finally, for all $j\in J$, one can follow the same argument as in the proof of Proposition~\ref{prop808} to attain bounds regarding $\left\|\cdot\right\|_{(1,R)}$. On the other hand, one has that
\begin{multline*}
\sup_{m\in\R,\tau\in S_d}\left|\mathcal{F}_j(m)\frac{\tau^j}{P_m(\tau)}\right|(1+|m|)^{\mu}e^{\beta|m|}\frac{1}{|\tau|}\exp\left(-\frac{k}{2\log(q)}\log^2|\tau|-\alpha\log|\tau|\right)\\
\le \left\|\mathcal{F}_j(m)\right\|_{(\beta,\mu)}\frac{1}{\delta_1 C_{R_{D}}}\max_{x>0}x^{j-1}\exp\left(-\frac{k}{2\log(q)}\log^2(x)-\alpha\log(x)\right)\le C_j\frac{H}{\delta_1 C_{R_{D}}},
\end{multline*}
where $C_{R_D}>0$ is such that $|R_{D}(im)|\le C_{R_D}(1+|m|)^{\hbox{deg}(R_D)}$ for all $m\in\R$, and with
$$H:=\max_{j\in J}\left(\max_{x>0}x^{j-1}\exp\left(-\frac{k}{2\log(q)}\log^2(x)-\alpha\log(x)\right)\right).$$

This entails that 
$$\left\|\mathcal{F}_j(m)\frac{\tau^j}{P_m(\tau)}\right\|_{(R,k,\alpha,\beta,\mu)}\le C_j\left(C_PR^j\sum_{p\ge0}\left(\frac{R}{R_1}\right)^{p}+\frac{H}{\delta_1 C_{R_D}}\right)\\
\le \varsigma_f\left(C_PR^j\sum_{p\ge0}\left(\frac{R}{R_1}\right)^{p}+\frac{H}{\delta_1 C_{R_D}}\right) .$$

Let $\varsigma_c,\varsigma_f>0$ be such that

\begin{multline*}
\sum_{\underline{\ell}=(\ell_0,\ell_1,\ell_2)\in\mathcal{A},\ell_2=1}\tilde{C}_3\frac{1}{(2\pi)^{1/2}q^{\frac{\ell_0(\ell_0-1)}{2k}}}+\sum_{\underline{\ell}=(\ell_0,\ell_1,\ell_2)\in\mathcal{A},\ell_2\ge2}\tilde{C}_4\frac{1}{(2\pi)^{1/2}}\frac{1}{q^{\frac{\ell_0(\ell_0-1)}{2k}}}+\\
\varsigma_f \left(C_P\sum_{p\ge0}\left(\frac{R}{R_1}\right)^{p}\left(\sum_{j\in J}R^j\right)+(\# J)\frac{H}{\delta_1 C_{R_D}}\right)\le 1.
\end{multline*}

Then, one has that $\mathcal{H}_2(\omega(\tau,m))\in E^{d}_{(R,k,\alpha,\beta,\mu)}$ with $\left\|\mathcal{H}_2(\omega(\tau,m))\right\|_{(R,k,\alpha,\beta,\mu)}\le \varpi$.

We now study the difference of the image of two elements by $\mathcal{H}_2$. More precisely, let $\omega_1(\tau,m),\omega_2(\tau,m)\in E^{d}_{(R,k,\alpha,\beta,\mu)}$ with $\left\|\omega_j(\tau,m)\right\|_{(R,k,\alpha,\beta,\mu)}\le \varpi$. 

From the proof of Proposition~\ref{prop808} we arrive at 
$$\left\|\mathcal{H}_2(\omega_1(\tau,m))-\mathcal{H}_2(\omega_2(\tau,m))\right\|_{(1,R)}\le \frac{1}{2}\left\|\omega_1(\tau,m)-\omega_2(\tau,m)\right\|_{(1,R)}.$$

On the other hand, analogous estimates as before yield

\begin{multline*}
\left\|\frac{1}{P_m(\tau)}\int_{-\infty}^{\infty}A_{\underline{\ell}}(m-m_1)\tau^{\ell_0}\sigma_{q;\tau}^{\ell_1-\frac{\ell_0}{k}}(\omega_1(\tau,m_1)-\omega_2(\tau,m_1))R_{\underline{\ell}}(im_1)dm_1\right\|_{(R,k,\alpha,\beta,\mu)}\\
\le \tilde{C}_3\left\|\omega_1(\tau,m)-\omega_2(\tau,m)\right\|_{(R,k,\alpha,\beta,\mu)},
\end{multline*}
together with
\begin{multline*}
\left\|\frac{1}{P_m(\tau)}\int_{-\infty}^{\infty}A_{\underline{\ell}}(m-m_1)\int_{\tilde{C}_{\rho/2}}\mathbb{D}_{q,k'}(x/\tau^{\ell_2})\sigma_{q;x}^{-k''}\left[x^{\ell_0}\sigma_{q;x}^{\ell_1-\frac{\ell_0}{k}}(\omega_1(x,m_1)-\omega_2(x,m_1))\right]\frac{dx}{x}\right.\\
\left.\times R_{\underline{\ell}}(im_1)dm_1\right\|_{(R,k,\alpha,\beta,\mu)}
\le \tilde{C}_4 \left\|\omega_1(\tau,m)-\omega_2(\tau,m)\right\|_{(R,k,\alpha,\beta,\mu)}.
\end{multline*}

We choose $\varsigma_c>0$ be such that

$$\sum_{\underline{\ell}=(\ell_0,\ell_1,\ell_2)\in\mathcal{A},\ell_2=1}\tilde{C}_3\frac{1}{(2\pi)^{1/2}q^{\frac{\ell_0(\ell_0-1)}{2k}}}+\sum_{\underline{\ell}=(\ell_0,\ell_1,\ell_2)\in\mathcal{A},\ell_2\ge2}\tilde{C}_4\frac{1}{(2\pi)^{1/2}}\frac{1}{q^{\frac{\ell_0(\ell_0-1)}{2k}}}\le \frac{1}{2},$$
leading to 
$$\left\|\mathcal{H}_2(\omega_1(\tau,m))-\mathcal{H}_2(\omega_2(\tau,m))\right\|_{(R,k,\alpha,\beta,\mu)}\le \frac{1}{2}\left\|\omega_1(\tau,m)-\omega_2(\tau,m)\right\|_{(R,k,\alpha,\beta,\mu)}.$$

The classical contractive mapping theorem applied to 
$$\mathcal{H}_2:\overline{D}(0,\varpi)\to \overline{D}(0,\varpi),$$
with $\overline{D}(0,\varpi)\subseteq E^{d}_{(R,k,\alpha,\beta,\mu)}$ allows us to conclude the result.
\end{proof}

\section{Main results}\label{secmainres}

In this section, we maintain the elements stated in Section~\ref{secmainpr} regarding the construction of the main problem under study, namely (\ref{epral}) under null initial conditions $u(0,z)\equiv 0$. In particular, we consider the geometric configuration of the problem described in Section~\ref{secgeo}.

We are in position to state the first main result of the present work.

\begin{theo}\label{teo1}
Let $\rho>0$ and let $S_d$ be an unbounded sector with vertex at the origin, and bisecting direction $d\in\R$, chosen in accordance to the geometric configuration described in Section~\ref{secgeo}. There exist $\varpi_{c}>0$ and $\varpi_{f}>0$ such that if 
$$C_{\underline{\ell}}\le \varpi_c\quad \underline{\ell}\in\mathcal{A},\qquad C_j\le \varpi_f,\quad j\in J,$$
($C_{\ell}$ being the constant determined in (\ref{e588}) and $C_{j}$ being the constant determined in (\ref{e596}), respectively)
then the problem (\ref{epral}) under null initial data $u(0,z)\equiv 0$ admits a formal solution  $\hat{u}(t,z)=\sum_{p\ge1}u_p(z)t^p\in (\mathcal{O}_b(H_{\beta'}))[[t]]$ which is $G_q-$summable of order $k$ with respect to $t$ in direction $d$.
\end{theo}
\begin{proof}
In the situation declared, Proposition~\ref{prop808} holds, which guarantees there exists $R\in (0,\rho)$ such that the operator $\mathcal{H}_1$ defined in (\ref{e614}) determines a fixed point $\omega(\tau,m)\in HE_{(\beta,\mu,R)}$. Therefore, according to Lemma~\ref{lema0}, $\omega(\tau,m)$ can be written in the form 
\begin{equation}\label{e969}
\omega(\tau,m)=\sum_{p\ge 1}\omega_p(m)\tau^p\in E_{(\beta,\mu)}\{\tau\},
\end{equation}
with $\omega_p(m)\in E_{(\beta,\mu)}$ and where the previous series converges on some neighborhood of the origin. Therefore, there exists $C_{\varsigma}>0$ such that
$$\left|\omega_p(m)\right|\le C_{\varsigma}\left(\frac{1}{R}\right)^{p}\frac{1}{(1+|m|)^{\mu}}e^{-\beta|m|},$$
for all $p\ge0$ and all $m\in\R$. We define the formal power series
\begin{equation}\label{e973}
\hat{U}(t,m)=\sum_{p\ge1}\omega_p(m)q^{\frac{p(p-1)}{2k}}t^p.
\end{equation}
It is clear from Corollary~\ref{coro467} that $\omega(\tau,m)$ stands for the $q$-Borel transform of order $k$ of (\ref{e973}), providing a formal solution to (\ref{eaux2}). In view of Proposition~\ref{prop882}, the definition of the function $\omega(\tau,m)$ can be prolonged to an infinite sector of bisecting direction $d$, with $q$-exponential growth at infinity. More precisely, for every $m\in\R$, the function $\tau\mapsto\omega(\tau,m)$ extends to a function defined in $S_d$. Moreover, this extension belongs to $E^{d}_{(R,k,\alpha,\beta,\mu)}$ for some small $R>0$. Therefore, if we maintain the same notation for this extension, there exists $C_{\varsigma,2}>0$ such that
\begin{equation}\label{e979}
|\omega(\tau,m)|\le C_{\varsigma,2}\frac{1}{(1+|m|)^{\mu}}e^{-\beta|m|}|\tau|\exp\left(\frac{k}{2\log(q)}\log^2|\tau|+\alpha\log|\tau|\right),
\end{equation}
for all $\tau\in S_d$ and all $m\in\R$. From the previous bounds, we deduce that the formal power series $\hat{U}(t,m)$ is $G_q$-summable of order $k$ with respect to $t$ in the direction $d$, with coefficients in the Banach space $E_{(\beta,\mu)}$ (see Section~\ref{sec56}). In addition to this, $\hat{U}(t,m)$ is a formal solution to (\ref{e169}), due to the computations made in Section~\ref{secI}.

Let us write $\hat{U}(t,m)=\sum_{p\ge1}U_p(m)t^p$. We define the formal power series
$$\hat{u}(t,z)=\sum_{p\ge1}u_p(z)t^p,$$
with 
$$u_p(z)=\frac{1}{(2\pi)^{1/2}}\int_{-\infty}^{\infty}U_p(m)e^{izm}dm,\quad p\ge1,$$
for all $z\in H_{\beta'}$, for any fixed $0<\beta'<\beta$. The strategy followed in Section~\ref{secI} guarantees that $\hat{u}(t,z)$ stands for a formal solution of (\ref{epral}). In view of the properties of $u_p(z)$ we derive that
$$|u_{p}(z)|\le \frac{1}{(2\pi)^{1/2}}C_{\varsigma}q^{\frac{p(p-1)}{2k}}\left(\frac{1}{R}\right)^{p}\int_{-\infty}^{\infty}\frac{1}{(1+|m|)^{\mu}}e^{-(\beta-\beta')|m|}dm,$$
for all $p\ge1$ and every $z\in H_{\beta'}$. This entails that $u_p(z)$ belongs to $(\mathcal{O}_{b}(H_{\beta'}),\left\|\cdot\right\|_{\infty})$ for all $p\ge1$. The series
$$\sum_{p\ge 1}\sup_{z\in H_{\beta'}}|u_p(z)|\frac{1}{q^{\frac{p(p-1)}{2k}}}(R')^{p}$$
is convergent for all $0<R'<R$. Therefore, the $q$-Borel transform of order $k$ along direction $d$ of $\hat{u}(t,z)$, 
$$\hat{\mathcal{B}}_{q;k}(t\mapsto\hat{u}(t,z))(\tau)=\sum_{p\ge1}\frac{u_p(z)}{q^{\frac{p(p-1)}{2k}}}\tau^p\in(\mathcal{O}_b(H_{\beta'}))[[\tau]]$$
is convergent in the disc $D(0,R')$, i.e. $\hat{\mathcal{B}}_{q;k}(t\mapsto\hat{u}(t,z))(\tau)\in\tau\mathcal{O}_b(H_{\beta'})\{\tau\}$. In accordance to (\ref{e969}), one has
$$\hat{\mathcal{B}}_{q;k}(t\mapsto\hat{u}(t,z))(\tau)=\frac{1}{(2\pi)^{1/2}}\int_{-\infty}^{\infty}\omega(\tau,m)e^{izm}dm,\quad z\in D(0,R'),$$
for every $0<R'<R$. In view of (\ref{e979}), one has
\begin{multline*}
\sup_{z\in H_{\beta'}}|\hat{\mathcal{B}}_{q;k}(t\mapsto\hat{u}(t,z))(\tau)|\le\frac{C_{\varsigma,2}}{(2\pi)^{1/2}}|\tau|\exp\left(\frac{k}{2\log(q)}\log^2|\tau|+\alpha\log|\tau|\right)\\
\times \left(\int_{-\infty}^{\infty}\frac{1}{(1+|m|)^{\mu}}e^{-(\beta-\beta')|m|}dm\right),
\end{multline*}
for $\tau\in S_d$. The fact that $\hat{\mathcal{B}}_{q;k}(t\mapsto\hat{u}(t,z))(\tau)\in\tau\mathcal{O}_b(H_{\beta'})\{\tau\}$ and the previous upper estimate in $S_d$ yields that the formal solution $\hat{u}(t,z)$ of (\ref{epral}) with $\hat{u}(0,z)\equiv 0$, is $G_q-$summable of order $k$ with respect to $t$ in direction $d$.
\end{proof}

The second main result of the present work describes a pseudo $q$-difference equation satisfied by the $G_q$-sum of order $k$ with respect to $t$ of the formal solution $\hat{u}(t,z)$ of (\ref{epral}), with $\hat{u}(t,0)\equiv 0$.

For that purpose, we define some operators acting on Fourier and $q$-Laplace transforms.

\begin{prop}\label{prop1006}
Let the hypothesis of Theorem~\ref{teo1} hold. Let $0<R<\rho$. For every $\boldsymbol{\omega}^d(\tau,m)\in E^d_{(R,k,\alpha,\beta,\mu)}$, let 
\begin{equation}\label{e1011}
\boldsymbol{u}^{d}(t,z)=\frac{\pi_{q,k}}{(2\pi)^{1/2}}\int_{-\infty}^{\infty}\int_{L_{d}}\Theta_{k}\left(\frac{t}{u}\right)\boldsymbol{\omega}^d(u,m)e^{imz}\frac{du}{u}dm,
\end{equation}
for $L_{d}=[0,\infty)e^{id}$, and where $\pi_{q,k}$ and $\Theta_{k}$ are defined in Proposition~\ref{prop191}.

Let us define the following operator acting on $\boldsymbol{u}^d(t,z)$:
\begin{equation}\label{e1017}
\hbox{Exp}_q\left(-\alpha_Dt^{d_D}\sigma_{q;t}^{\frac{d_D}{k}}\right)\boldsymbol{u}^{d}(t,z):=\frac{\pi_{q,k}}{(2\pi)^{1/2}}\int_{-\infty}^{\infty}\int_{L_{d}}\frac{1}{\exp_q(\tilde{\alpha}_Du^{d_D})}\Theta_{k}\left(\frac{t}{u}\right)\boldsymbol{\omega}^d(u,m)e^{imz}\frac{du}{u}dm,
\end{equation}
where $\tilde{\alpha}_D>0$ is given by (\ref{e382}).

For every $\underline{\ell}=(\ell_0,\ell_1,\ell_2)\in\mathcal{A}$ with $\ell_2\ge2$ we define the next operator acting on $\boldsymbol{u}^d(t,z)$:
\begin{multline}\label{e1018}
\mathbb{G}_{\underline{\ell},k}(\boldsymbol{u}^{d})(t,z):=\frac{\pi_{q,k}}{(2\pi)^{1/2}}\int_{-\infty}^{\infty}\int_{L_{d}}\frac{1}{\exp_q(\tilde{\alpha}_D u^{d_D})}\Xi_{q,k'}\int_{\tilde{C}_{\rho/2}}\mathbb{D}_{q,k'}(x/ u^{\ell_2})\\
\times\sigma_{q;x}^{-k''}\left[x^{\ell_0}\sigma_{q;x}^{\ell_1-\frac{\ell_0}{k}}\boldsymbol{\omega}(x,m)\right]\frac{dx}{x} \Theta_{k}\left(\frac{t}{u}\right)e^{imz}\frac{du}{u}dm,
\end{multline}
provided that $d_D$ satisfies the inequality (\ref{e676c}).

Then, the functions $\hbox{Exp}_q(-\alpha_D t^{d_D}\sigma_{q;t}^{\frac{d_D}{k}})\boldsymbol{u}^d(t,z)$ and $\mathbb{G}_{\underline{\ell},k}(\boldsymbol{u}^{d})(t,z)$ are well defined and are bounded holomorphic on $\tilde{D}(0,R_1)\times H_{\beta'}$ for some radius $R_1>0$ and any fixed $0<\beta'<\beta$.
\end{prop}
\begin{proof}
Observe that (\ref{e1011}) is well-defined due to $\boldsymbol{\omega}^d(\tau,m)\in E^d_{(R,k,\alpha,\beta,\mu)}$, and $\boldsymbol{u}^{d}(t,z)$ turns out to be a holomorphic and bounded function on $\tilde{D}(0,R_1)\times H_{\beta'}$, for some $R_1>0$ and any fixed $0<\beta'<\beta$, in accordance to Proposition~\ref{prop191} and Section~\ref{ap2}.

It is also worth remarking that the operator (\ref{e1017}) is well defined applied on functions $\boldsymbol{u}^d$ of the form (\ref{e1011}). It is also important to point out that the function $\tau\mapsto 1/\exp_q(\tilde{\alpha}_D\tau^{d_D})$ exists along direction $L_{d}$ taking into account that $\exp_q$ is an entire function with $\exp_q(0)=1$ whose zeros are located on the negative real axis. Moreover, it admits a subexponential growth at infinity (see (\ref{e411}) and Proposition~\ref{prop430}). Analogous estimates as the one of (\ref{e422}) in Lemma~\ref{lema486b} guarantee that the operator $\mathbb{G}_{\underline{\ell},k}$ is well-defined when acting on map $\boldsymbol{u}^d$ of the form (\ref{e1011}). In both situations, the result of the operator is a holomorphic and bounded function on $\tilde{D}(0,R_1)\times H_{\beta'}$, for some $R_1>0$ and any fixed $0<\beta'<\beta$.
\end{proof}

We are in conditions to state the second main result of the present work.

\begin{theo}\label{teo2}
Let the hypothesis of Theorem~\ref{teo1} hold. We consider the formal solution $\hat{u}(t,z)$ of (\ref{epral}) with null initial data $\hat{u}(0,z)\equiv 0$. Theorem~\ref{teo1} guarantees that $\hat{u}(t,z)$ is $G_q-$summable of order $k$ with respect to $t$ in direction $d$. Let $u^d(t,z)$ its $G_q-$sum along direction $d$, which is a holomorphic and bounded map defined on $\tilde{D}(0,R_1)\times H_{\beta'}$, for some $R_1>0$ and any fixed $0<\beta'<\beta$.

Then, $u^d(t,z)$ satisfies the following pseudo $q$-difference-differential equation:
\begin{multline}\label{e1036}
\hbox{Exp}_q\left(-\alpha_Dt^{d_{D}}\sigma_{q;t}^{\frac{d_D}{k}}\right)Q(\partial_z)u^d(t,z)=R_D(\partial_z)u^d(t,z)\\
+\sum_{\underline{\ell}=(\ell_0,\ell_1,\ell_2)\in\mathcal{A},\ell_2=1}a_{\underline{\ell}}(z)R_{\underline{\ell}}(\partial_z)\hbox{Exp}_q\left(-\alpha_Dt^{d_{D}}\sigma_{q;t}^{\frac{d_D}{k}}\right)\left(t^{\ell_0}\sigma_{q;t}^{\ell_1}u^d\right)(t,z)\\
+\sum_{\underline{\ell}=(\ell_0,\ell_1,\ell_2)\in\mathcal{A},\ell_2\ge2}a_{\underline{\ell}}(z)R_{\underline{\ell}}(\partial_z)\mathbb{G}_{\underline{\ell},k}(u^{d})(t,z)+\hbox{Exp}_q\left(-\alpha_Dt^{d_{D}}\sigma_{q;t}^{\frac{d_D}{k}}\right)f(t,z),
\end{multline}
on $\tilde{D}(0,R_1)\times H_{\beta'}$ for small enough $R_1>0$ and any $0<\beta'<\beta$.
\end{theo} 
\begin{proof}

Taking into account the assumptions made in Theorem~\ref{teo1} one can apply Proposition~\ref{prop1006}. Departing from the solution $\omega^d(\tau,m)$ of the auxiliary equation (\ref{eaux2}), which belongs to the Banach space $E^{d}_{(R,k,\alpha,\beta,\mu)}$ with $\left\|\omega^d(\tau,m)\right\|_{(R,k,\alpha,\beta,\mu)}\le \varpi$ for some fixed $\varpi>0$.

The $G_q$-sum of order $k$ with respect to $t$ of the formal solution, $\hat{u}(t,z)$ to (\ref{epral}) can be written in the form of an inverse Fourier and $q$-Laplace transform (\ref{e1011}), defining a holomorphic and bounded function on $\tilde{D}(0,R)\times H_{\beta'}$. 

After multiplying the auxiliary equation (\ref{eaux2}) satisfied by $\omega^d(\tau,m)$ by $P_m(\tau)/\exp_{q}(\tilde{\alpha}_d\tau^{d_D})$ one arrives at the equation

\begin{multline}\label{eaux2b}
\frac{Q(im)}{\exp_{q}(\tilde{\alpha}_D\tau^{d_D})}\omega^d(\tau,m)=R_D(im)\omega^d(\tau,m)\\
+\sum_{\underline{\ell}=(\ell_0,\ell_1,\ell_2)\in\mathcal{A},\ell_2=1}\frac{1}{(2\pi)^{1/2}}\int_{-\infty}^{\infty}\frac{A_{\underline{\ell}}(m-m_1)}{\exp_{q}(\tilde{\alpha}_D\tau^{d_D})}\frac{\tau^{\ell_0}}{q^{\frac{\ell_0(\ell_0-1)}{2k}}}\sigma_{q;\tau}^{\ell_1-\frac{\ell_0}{k}}\omega^d(\tau,m_1)R_{\underline{\ell}}(im_1)dm_1\\
+\sum_{\underline{\ell}=(\ell_0,\ell_1,\ell_2)\in\mathcal{A},\ell_2\ge2}\frac{1}{(2\pi)^{1/2}}\int_{-\infty}^{\infty}\frac{A_{\underline{\ell}}(m-m_1)}{\exp_{q}(\tilde{\alpha}_D\tau^{d_D})}\Xi_{q,k'}\int_{\tilde{C}_{\rho/2}}\mathbb{D}_{q,k'}(x/\tau^{\ell_2})\\
\times\sigma_{q;x}^{-k''}\left[\frac{x^{\ell_0}}{q^{\frac{\ell_0(\ell_0-1)}{2k}}}\sigma_{q;x}^{\ell_1-\frac{\ell_0}{k}}\omega^d(x,m_1)\right]\frac{dx}{x} R_{\underline{\ell}}(im_1)dm_1+\frac{1}{\exp_{q}(\tilde{\alpha}_D\tau^{d_D})}\sum_{j\in J}\mathcal{F}_j(m)\tau^j,
\end{multline}
valid for $\tau\in S_d\cup D(0,R)$ for some $0<R<\rho$, and $m\in\R$. We apply $q$-Laplace transform of order $k$ along direction $d$ and inverse Fourier transform to the previous equation, to arrive at equation (\ref{e1036} being satisfied by $u^{d}(t,z)$ in $\tilde{D}(0,R_1)\times H_{\beta'}$ for small enough $R_1>0$. At this point we have made use of the notation of the operators introduced in Proposition~\ref{prop1006}.
\end{proof}

\vspace{0.3cm}

\textbf{Remark:} Observe that the map $u^{d}(t,z)$ does not in general satisfy equation (\ref{epral}), but (\ref{e1036}):
\begin{itemize}
\item The infinite order $q$-difference operator $\exp_q\left(\alpha_Dt^{d_{D}}\sigma_{q;t}^{\frac{d_D}{k}}\right)$ can not be directly applied to $u^{d}(t,z)$ since the map $\tau\mapsto \exp(\alpha_D\tau^{d_D})$ has $q$-exponential growth at infinity of order $d_D^2$ (see (\ref{e398b}) and its consequence) whereas the limit growth rate at infinity admitted for a function in $S_d$ is smaller, as described in Section~\ref{sec56}.  
\item Mahler operator $t\mapsto t^{\ell_2}$ for all $\underline{\ell}=(\ell_0,\ell_1,\ell_2)\in\mathcal{A}$ can not be applied to $t^{\ell_0}\sigma_{q;t}^{\ell_1}u^{d}(t,z)$ as it is associated to the mapping 
$$\tau\mapsto \int_{\tilde{C}_{\rho/2}}\mathbb{D}_{q,k'}(x/\tau^{\ell_2})\sigma_{q;x}^{-k''}\left[x^{\ell_0}\sigma_{q;x}^{\ell_1-\frac{\ell_0}{k}}\omega^d(x,m_1)\right]\frac{dx}{x},$$
whose growth at infinity is governed by a $q$-exponential growth of order $\frac{k\ell_2^2}{\ell_2^2-1}$, owing to the bound (\ref{e422}), exceeding the limit $q$-exponential growth rate of order $k$.
\end{itemize}

\section{Appendix I: Review on $q$-analogs}\label{secapqex}

This appendix aims to collect some of the main facts on the $q$-analogs of different analytic tools which appear in the present work. More precisely, the $q$-exponential function, and a continuous $q$-analog of Laplace transform of some positive order. In the whole section, we fix $q>1$ and $k>0$.

\subsection{$q$-exponential function}\label{secapqex1}

This section is devoted to recall the definition and main facts on the $q$-exponential function. We refer to~\cite{tahara} and the references therein for further details and related topics.

A $q$-analog of the exponential function is defined by
\begin{equation}\label{e574}
\exp_q(z)=\sum_{n\ge0}\frac{z^n}{[n]_q^{!}},
\end{equation}
where $[n]_{q}^{!}$ stands for the $q$-factorial of the nonnegative integer $n$, given by
$$[0]_q^{!}=1,\qquad [n]_{q}^{!}=[1]_q[2]_q\cdots [n]_q.$$
Here, $[j]_q$ stands for the $q$-number defined by $[j]_q=1+q+\ldots+q^{j-1}$ for every positive integer $j$.

The function $\exp_q(z)$ turns out to be an entire function on $\C$, which has simple zeros at $z=-q^{m+1}/(q-1)$, for every nonnegative integer $m$. Therefore, we observe there exists $C_0>0$ such that 
$$|\exp_q(z)|>C_0,\quad z\in D(0,\frac{q^{1/2}}{q-1}).$$ 
The location of zeros also motivates the definition of the following sets, isolating neighborhoods of the zeros of $\exp_q$. Given $\epsilon>0$ and a nonnegative integer $m$. We denote 
$$B_{m,\epsilon}=\left\{z\in\C: \left|z+\frac{q^{m+1}}{q-1}\right|\le \epsilon \frac{q^{m+1}}{q-1}\right\},$$
and $\mathcal{L}_{\epsilon}=\cup_{m\ge0}B_{m,\epsilon}$.

In accordance to Proposition 2.2 (4)~\cite{tahara}, one has that given $\epsilon>0$, the map $z\mapsto \exp_q(z)$ is rapidly increasing in $\C\setminus\mathcal{L}_{\epsilon}$ as $|z|\to+\infty$. More precisely, there exists $K_1>0$ such that 
\begin{equation}\label{e411}
|\exp_q(z)|\le K_1\exp(\mu(|z|))=K_1\exp\left(\frac{\log^2|z|}{2\log(q)}+\left(-\frac{1}{2}+\frac{\log(q-1)}{\log(q)}\right)\log|z|\right),
\end{equation}
for all $z\in\C\setminus\mathcal{L}_\epsilon$ with $|z|\ge \frac{q^{1/2}}{q-1}$. In addition to this, there exists $K_0>0$ such that
$$|\exp_q(z)|\ge \frac{\epsilon}{K_0}\exp(\mu(|z|)),$$
for all $z\in\C\setminus\mathcal{L}_\epsilon$ with $|z|\ge \frac{q^{1/2}}{q-1}$.

A geometric consequence of the previous properties is the following. Let $0<\epsilon<1$ and $m\ge0$ be an integer. Let $\theta\in (0,\pi/2)$ such that $\sin(\theta)=\epsilon$. Given $d(m)=\epsilon\frac{q^{m+1}}{q-1}$ one observes that the sector 
\begin{equation}\label{e419}
S_{0,\pi-\theta}=\{z\in\C:|\hbox{arg}(z)|<\pi-\theta\}
\end{equation}
and the disc $D(-\frac{q^{m+1}}{q-1},d(m))$ have empty intersection. As a result, for small enough $\theta>0$, one can find small enough $\epsilon>0$ such that 
$$S_{0,\pi-\theta}\subseteq \C\setminus\mathcal{L}_{\epsilon}.$$

\begin{figure}
	\centering
\includegraphics[width=0.45\textwidth]{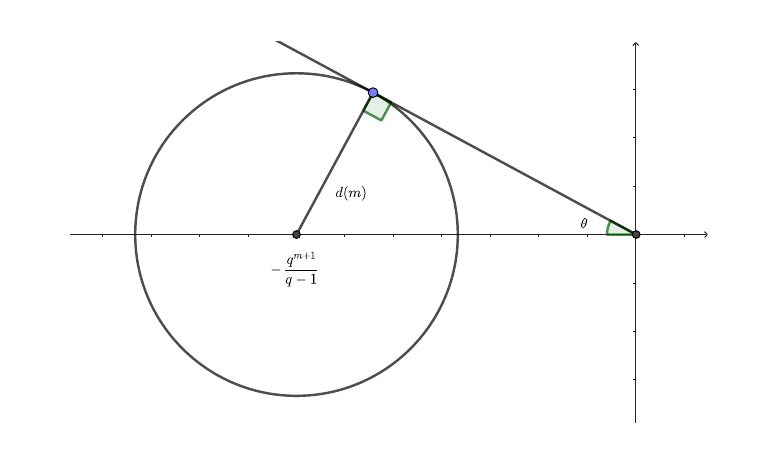}
\caption{Geometric distributions of zeros of $\exp_q(z)$}
\label{fig1}
\end{figure}

As a consequence, the following result holds.

\begin{prop}\label{prop430}
Let $\theta\in(0,\pi/2)$ be a given angle and let $\epsilon>0$ such that $\sin(\theta)=\epsilon$. Then, for any given unbounded sector $U$ with vertex at the origin with $U\subseteq S_{0,\pi-\theta}$, there exist $K_0,K_1>0$  and $C_0>0$ such that the following estimates hold.
\begin{itemize}
\item $|\exp_q(z)|\le K_1\exp(\mu(|z|))$, for all $z\in U\setminus D(0,\frac{q^{1/2}}{q-1})$,
\item $|\exp_q(z)|\ge \frac{\epsilon}{K_0}\exp(\mu(|z|))$, for all $z\in U\setminus D(0,\frac{q^{1/2}}{q-1})$,
\item $|\exp_q(z)|\ge C_0$, for all $z\in \overline{D}(0,\frac{q^{1/2}}{q-1})$,
\end{itemize}
with $\mu(\cdot)$ defined by (\ref{e411}).
\end{prop}

\subsection{On $q$-Laplace transform}\label{sec200}

One can find different definitions in the literature of $q$-analogs of Laplace transform in both discrete and continuous form. We refer to~\cite{vizh} for the definition and main properties of several of them. The definition and main properties of the $q$-analog of the Laplace transform of positive order $k$ considered in this work was introduced in~\cite{zh99,mazh}. We include the proofs of the main results for the sake of completeness.

\begin{prop}\label{prop191}
Let $\rho>0$, and $S_d$ be an unbounded sector with vertex at the origin and bisecting direction $d\in\R$. Let $f:D(0,\rho)\cup S_d\to\C$ be a holomorphic function on $D(0,\rho)\cup S_d$, continuous on $\overline{D}(0,\rho)$. We assume there exist $K,\alpha>0$ and $\delta>1$ such that
\begin{equation}\label{e193}
|f(z)|\le K|z|\exp\left(\frac{k}{2}\frac{\log^2(|z|+\delta)}{\log(q)}+\alpha\log(|z|+\delta)\right),
\end{equation}
for all $z\in D(0,\rho)\cup S_d$. Let $\theta_d$, with $\{re^{i\theta_d}:r>0\}\subseteq S_{d}$. 

The $q$-Laplace transform of order $k$ of $f$ along direction $\theta_d$ is defined as
$$\mathcal{L}^{\theta_d}_{q;k}(f(z))(T)=\pi_{q,k}\int_{L_{\theta_d}}\Theta_k\left(\frac{T}{u}\right)f(u)\frac{du}{u}$$
where
$$\pi_{q,k}:=\frac{q^{-\frac{1}{8k}}\sqrt{k}}{\sqrt{2\pi \log(q)}},$$
the kernel function is defined by  
$$\Theta_k(z):=q^{-\frac{k}{2}\frac{\log(z)}{\log(q)}\left(\frac{\log(z)}{\log(q)}-\frac{1}{k}\right)},$$
and the integration path is given by $L_{\theta_d}=[0,\infty)e^{i\theta_d}$.

In the previous situation, the function $T\mapsto \mathcal{L}^{\theta_d}_{q;k}(f(z))(T)$ is analytic on the universal covering $\tilde{D}(0,R)$ of the disc $D(0,R)$, for some small enough $R>0$. 

In addition to this, $\mathcal{L}^{\theta_d}_{q;k}(f(z))(T)$ is bounded on every sectorial domain $\tilde{S}_{I,R}=\{re^{i\theta}\in\mathcal{R}:0<r<R,\theta\in I\}$, for some bounded segment $I\subseteq\R$. More precisely, there exists $K>0$ such that for all $T\in \tilde{D}(0,R)$ and for all $\theta_d\in(d-\delta_1,d+\delta)$ (with $e^{i\theta_d}\in S_d$), then
$$\left|\mathcal{L}^{\theta_d}_{q;k}(f(z))(T)\right|\le Ke^{\frac{k}{2\log(q)}(\hbox{arg}(T)-\theta_d)}.$$
\end{prop}

\textbf{Remark:} It is worth mentioning that the value of $\mathcal{L}^{\theta_d}_{q;k}(f(z))(T)$ does not depend on the choifce of $\theta_d$ provided that $L_{\theta_d}\subseteq S_d\cup\{0\}$, as it follows from the holomorphy of $f$ and Cauchy's theorem.

\begin{proof}
The proof follows from classical computations. We only sketch some of its points for the case $d=0$, whereas other choice of $d\in\R$ can be handled analogously. We first observe that the expression 
$$\left|e^{-\frac{k}{2\log(q)}\log^2(T/u)}\right|=e^{\frac{k}{2\log(q)}\log(u)\log|T|}e^{-\frac{k}{2\log(q)}\log^2(u)}$$
and
$$\left|e^{\frac{1}{2}\log(T/u)}\right|=e^{\frac{1}{2}\log(|T|/u)}$$
for any $u>0$, and $T\in\C\setminus\{0\}$.

Two cases arise:

\begin{itemize}
\item[$-$] Case 1: take $0\le u\le 1$. By construction, there exists $K>0$ such that $|f(u)/u|\le K$ for all $0\le u\le 1$. We observe that 
$$e^{\frac{1}{2}\log(|T|/u)}e^{-\frac{k}{2\log(q)}\log^2(\frac{|T|}{u})}e^{\frac{k}{2\log(q)}(\hbox{arg}(T))^2}$$
is bounded by some constant $K_1>0$ provided that $\hbox{arg}(T)\in I$ for any bounded interval $I\subseteq \R$ and all $T\in\C$ with $|T|\le 1$. As a result,  there exists $K_2>0$ such that 
\begin{equation}\label{e1204a}
\left|e^{-\frac{k}{2\log(q)}\log^2(T/u)+\frac{1}{2}\log(T/u)}\right|\left|\frac{f(u)}{u}\right|\le K_2,
\end{equation}
for all $0\le u\le 1$ and $|T|\le 1$.

\item[$-$] Case 2: take $u\ge 1$. Since $\log(1+x)\equiv x$ for $x\to 0$ and $\lim_{x\to\infty}\log(u)/u=0$, one gets the existence of $\tilde{K}_1,\tilde{K}_2>0$ such that
$$\log^2(u+\delta)\le \log^2(u)+\tilde{K}_1,\quad \log(u+\delta)\le \log(u)+\tilde{K}_2,$$
for $u\ge 1$. As a result,
$$\left|\frac{f(u)}{u}\right|\le K\exp\left(\frac{k}{2\log(q)}\log^2(u)+\frac{\tilde{K}_1k}{2\log(q)}+\alpha\log(u)+\alpha\tilde{K}_2\right)$$
for $u\ge1$. We assume that $|T|\le R<1$. Then, one has that
$$e^{\frac{k}{\log(q)}\log|T|\log(u)}\le e^{\frac{k}{\log(q)}\log(R)\log(u)}.$$
As a result,
\begin{multline}\label{e1204b}
\left|e^{-\frac{k}{2\log(q)}\log^2(T/u)+\frac{1}{2}\log(T/u)}\right|\left|\frac{f(u)}{u}\right|\\
\le e^{-\frac{k}{2\log(q)}\log^2|T|}e^{\frac{k}{2\log(q)}(\hbox{arg}(T))^2}e^{\frac{k}{2\log(q)}\log(R)\log(u)}e^{-\frac{k}{2\log(q)}\log^2(u)}e^{\frac{1}{2}\log|T|}e^{-\frac{1}{2}\log(u)}\\
\times K\exp\left(\frac{k}{2\log(q)}\log^2(u)+\frac{\tilde{K}_1k}{2\log(q)}+\alpha\log(u)+\alpha\tilde{K}_2\right).
\end{multline}
Observe that if one takes $0<R<1$ in such a way that 
$$\frac{k}{\log(q)}\log(R)-\frac{1}{2}+\alpha<-\frac{3}{2},$$
then 
$$e^{\frac{k}{\log(q)}\log(R)\log(u)-\frac{1}{2}\log(u)+\alpha\log(u)}\le e^{-\frac{3}{2}\log(u)}=\frac{1}{u^{3/2}}.$$
As a result, one gets that
$$\int_1^{+\infty}\left|e^{-\frac{k}{2\log(q)}\log^2(T/u)+\frac{1}{2}\log(T/u)}\right|\left|\frac{f(u)}{u}\right|du\le KM_1e^{\frac{k}{2\log(q)}(\hbox{arg}(T))^2}e^{\frac{\tilde{K}_1k}{2\log(q)}+\alpha \tilde{K}_2}\int_{1}^{+\infty}\frac{1}{u^{3/2}}du$$
for some $M_1>0$, valid for all $t\in \tilde{D}(0,R)$.
\end{itemize}
Taking into account (\ref{e1204a}) and (\ref{e1204b}) one concludes the proof.
\end{proof}

\subsection{On the action of $q$-difference operators and $q$-Laplace transform}

We describe the action of $q$-difference operators acting on the $q$-Laplace transform of order $k$ of some function.

\begin{prop}
Let $f:D(0,\rho)\cup S_d\to \C$ be a function satisfying the requirements of Proposition~\ref{prop191}. Let $\sigma\ge 0$, $j\ge0$ be integer numbers. Then, it holds that
$$ T^{\sigma}\sigma_{q,T}^{j}\left(\mathcal{L}^{d}_{q;k}(f(z))\right)(T)=\mathcal{L}^{d}_{q;k}\left(\frac{z^{\sigma}}{(q^{1/k})^{\sigma(\sigma-1)}}\sigma_{q,z}^{j-\frac{\sigma}{k}}f(z)\right)(T),$$
for every $T\in \tilde{D}(0,R)$ for some small enough $R>0$.
\end{prop}
\begin{proof}
By definition, one has that

$$\mathcal{L}^{d}_{q;k}(z^{\sigma}f(q^{j-\frac{\sigma}{k}}z)(T)=\pi_{q,k}\int_{L_d}e^{-\frac{k}{2}\frac{\log^2(T/u)}{\log(q)}+\frac{1}{2}\log(T/u)}u^{\sigma}f(q^{j-\frac{\sigma}{k}}u)\frac{du}{u}.$$
The change of variable $u'=q^{j- \frac{\sigma}{k}}u$ transforms the previous expression into
$$\pi_{q,k}\int_{L_d}e^{-\frac{k}{2}\frac{\log^2(\frac{T}{q^{\sigma/k-j}u'})}{\log(q)}+\frac{1}{2}\log(\frac{T}{q^{\sigma/k-j}u'})}q^{(\sigma/k-j)\sigma}(u')^{\sigma}f(u')\frac{du'}{u'}.$$

Usual computations allow us to write
$$\exp\left(-\frac{k}{2}\frac{\log^2(\frac{T}{q^{\sigma/k-j}u'})}{\log(q)}+\frac{1}{2}\log(\frac{T}{q^{\sigma/k-j}u'})\right)=\Theta_k\left(\frac{T}{u'}\right)H_1 H_2,$$
where
$$H_1=\exp\left(-\frac{k}{2\log(q)}\log^2(\frac{q^jT}{u'})+\frac{1}{2}\log(\frac{q^jT}{u'})\right),$$
and
$$H_2=\exp\left(-\frac{k}{\log(q)}\log(q^{-\frac{\sigma}{k}})\log(\frac{q^jT}{u'})\right).$$

We observe after usual simplifications that 
$$H_1=q^{-\frac{\sigma^2}{2k}}q^{-\frac{\sigma}{2k}},\quad H_2=q^{\sigma j}T^{\sigma}(u')^{-\sigma}.$$

As a result, one concludes that 
$$\mathcal{L}^{d}_{q;k}(z^{\sigma}f(q^{j-\frac{\sigma}{k}}z)(T)=q^{\frac{\sigma^2}{2k}-\frac{\sigma}{2k}}T^{\sigma}\sigma_{q;T}^{j}\left(\mathcal{L}^{d}_{q;k}(f(z))(T)\right).$$
\end{proof}

Concerning the action of $q$-Laplace transform on polynomials, one can state the following result.

\begin{prop}\label{prop254}
Let $n\ge0$ be an integer. Then, it holds that
$$\mathcal{L}_{q;k}^d(z^n)(T)=q^{\frac{n(n-1)}{2k}}T^n.$$
\end{prop}
\begin{proof}
Let us assume without loss of generality that $d=0$, and $T\in (0,\infty)$. The general case is a direct consequence of Cauchy's formula.

Let us denote 
$$L(T)=\int_0^{+\infty}\exp\left(-\frac{k}{2}\frac{\log^2(\frac{T}{u})}{\log(q)}+\frac{1}{2}\log\left(\frac{T}{u}\right)\right)u^n\frac{du}{u}.$$
The change of variable $u'=T/u$ leads us to
$$L(T)=T^n\int_0^{+\infty}\exp\left(-\frac{k}{2\log(q)}\log^2(u')\right)\frac{1}{(u')^{n+1/2}}du'.$$
A second change of variable $h=\left(\frac{k}{2\log(q)}\right)^{1/2}\log(u')$ transforms the previous expression of $L(T)$ into
$$T^n\left(\frac{2\log(q)}{k}\right)^{1/2}\int_{-\infty}^{+\infty}e^{-h^2}e^{-a(n)h}dh,$$
with $a(n)=(-n+1/2)(2\log(q)/k)^{1/2}$. The next Gaussian equality, see~\cite{ma11}, can be applied
$$e^{a^2/4}\sqrt{\pi}\int_{-\infty}^{\infty}e^{-x^2-ax}dx,$$
in order to achieve
$$L(T)=T^n\left(\frac{2\log(q)}{k}\right)^{1/2}\sqrt{\pi}\exp\left(\frac{\log(q)}{2k}(-n+\frac{1}{2})^{2}\right).$$
We conclude that
$$L(T)=T^n\left(\frac{2\log(q)}{k}\right)^{1/2}\sqrt{\pi}q^{\frac{1}{8k}}q^{\frac{n^2-n}{2k}},$$
which is an expression equivalent to the statement in the proposition.
\end{proof}

\subsection{Asymptotic expansions of the $q$-Laplace transform near the origin}

Let $\rho>0$ and $d\in\R$. Let us also assume that $f:D(0,\rho)\cup S_d\to\C$ has the following Taylor expansion at the origin
$$f(z)=\sum_{n\ge1}\frac{a_n}{q^{\frac{n(n-1)}{2k}}}z^n,$$
for certain $a_n\in\C$, $n\ge 1$. We moreover assume that $f$ satisfies the properties stated in Proposition~\ref{prop191}. Then, the following result holds.

\begin{prop}\label{prop325}
The function $\mathcal{L}_{q;k}^{d}(f(z))(T)$ admits the formal power series $\hat{\mathcal{L}}(T)=\sum_{n\ge1}a_nT^n$ as its $q$-Gevrey asymptotic expansion of order $k$ on any sectorial domain of the form $\tilde{S}_R=\{\rho e^{i\theta}:0<\rho < R, \theta\in I\}$, for some bounded segment $I\subseteq\R$. In other words, there exist $K,A>0$ such that
$$\left|\mathcal{L}_{q;k}^{d}(f(z))(T)-\sum_{n=1}^{N-1}a_nT^n\right|\le K A^N q^{\frac{N^2}{2k}}|T|^{N},$$
for every $N\ge 2$ and all $T\in \tilde{S}_R$. 
\end{prop}

\subsection{Analytic $q$-Borel transform of order $k$}\label{sec332}

In this section, we recall the definition and some of the main analytic properties of $q$-Borel transform of some positive order. We describe the application of this transformation to polynomials.
 
\begin{defin}
Let $T\mapsto \varphi(T)$ be an analytic function on $\tilde{D}(0,R)$, for some $R>0$. Moreover, we assume that for all $re^{it}\in \tilde{D}(0,R)$, $\varphi$ satisfies the following bounds
\begin{equation}\label{e415}
\left|\varphi(re^{i t})\right|\le K e^{\frac{k}{2\log(q)}(t-\theta)^2},
\end{equation}
for some $K>0$ which does not depend on $r,t$, valid for all $\theta \in I=[d-\delta_1,d+\delta_1]$ for some small positive $\delta$.

We define the analytic $q$-Borel transform of $\varphi$ of order $k$ as the map
$$\mathcal{B}_{q;k}^{d}(\varphi)(\xi)=-\frac{q^{\frac{1}{8k}}\sqrt{k}i}{\sqrt{2\pi \log(q)}}\int_{\tilde{C}_{R/2}}q^{\frac{k}{2}\frac{\log(x/\xi)}{\log(q)}\left(\frac{\log(x/\xi)}{\log(q)}-\frac{1}{k}\right)}\varphi(x)\frac{dx}{x},$$
where $\tilde{C}_{R/2}=\left\{\frac{R}{2}e^{i\theta}:\theta \in\R\right\}\subseteq\tilde{D}_{R}$.
\end{defin}

\begin{prop}\label{prop302}
The map $\xi\mapsto \mathcal{B}_{q;k}^d(\varphi)(\xi)$ defines an analytic function on the sectorial domain 
$$V_I=\{re^{i\theta}:r\ge0, \theta\in I\}.$$
Furthermore, there exist $K,\alpha>0$ and $\delta>1$ such that
$$\left|\mathcal{B}_{q;k}^d(\varphi)(\xi)\right|\le K\exp\left(\frac{k}{2\log(q)}\log^2(|\xi|+\delta)+\alpha\log(|\xi|+\delta)\right),$$
for all $\xi\in V_I$ with $|\xi|\ge 1$.
\end{prop}

\textbf{Remark:} It is worth remarking that the map $\xi\mapsto \mathcal{B}_{q;k}^{d}(\varphi)(\xi)$ is univalued whenever $\varphi$ is univalued. Indeed, observe that under this assumption, one has after writing $x=R/2e^{i t}$ and the change of variables $s=t-2\pi$, that
\begin{multline*}
\mathcal{B}_{q;k}^{d}(\varphi)(\xi e^{i 2\pi})=-\frac{q^{\frac{1}{8k}}\sqrt{k}i}{\sqrt{2\pi \log(q)}}\int_{\R}q^{\frac{k}{2}\frac{\log(R/(2\xi)e^{i(t-2\pi)})}{\log(q)}\left(\frac{\log(R/(2\xi)e^{i(t-2\pi)})}{\log(q)}-\frac{1}{k}\right)}\varphi(\frac{R}{2}e^{it})idt\\
=-\frac{q^{\frac{1}{8k}}\sqrt{k}i}{\sqrt{2\pi \log(q)}}\int_{\R}q^{\frac{k}{2}\frac{\log(R/(2\xi)e^{is})}{\log(q)}\left(\frac{\log(R/(2\xi)e^{is})}{\log(q)}-\frac{1}{k}\right)}\varphi(\frac{R}{2}e^{is}e^{i 2\pi})idt,
\end{multline*}
which concludes with the fact that $\varphi(\frac{R}{2}e^{is}e^{i 2\pi})=\varphi(\frac{R}{2}e^{is})$, and therefore
$$\mathcal{B}_{q;k}^{d}(\varphi)(\xi e^{i 2\pi})=\mathcal{B}_{q;k}^{d}(\varphi)(\xi).$$

\begin{proof}
We observe that for $x=\frac{R}{2}e^{it}$, $t\in\R$, one has that
$$\frac{x}{\xi}=\frac{R}{2|\xi|}e^{i(t-\theta)},$$
together with
$$\log^2(\frac{x}{\xi})=\log^2\left(\frac{R}{2|\xi|}\right)+2i(t-\theta)\log\left(\frac{R}{2|\xi|}\right)-(t-\theta)^2.$$
The previous equalities yield
\begin{multline*}
\left|q^{\frac{k}{2}\frac{\log(x/\xi)}{\log(q)}\left(\frac{\log(x/\xi)}{\log(q)}-\frac{1}{k}\right)}\right|=\left|\exp\left(\frac{k}{2\log(q)}\log^2(x/\xi)-\frac{1}{2}\log(x/\xi)\right)\right|\\
=\exp\left(\frac{k}{2\log(q)}\log^2\left(\frac{k}{2|\xi|}\right)-\frac{k}{2\log(q)}(t-\theta)^2-\frac{1}{2}\log\left(\frac{k}{2|\xi|}\right)\right).
\end{multline*}
Classical computations guarantee that, provided that $R>0$ is small enough, there exist positive $K_1,\alpha_1$ such that
\begin{multline*}
\left|\exp\left(\frac{k}{2\log(q)}\log^2\left(\frac{x}{\xi}-\frac{1}{2}\log\left(\frac{x}{\xi}\right)\right)\right)\right|\\
\le K_1\exp\left(\frac{k}{2\log(q)}\log^2|\xi|\right)\exp\left(\alpha_1\log|\xi|\right)\exp\left(-\frac{k}{2\log(q)}(t-\theta)^2\right),
\end{multline*}
for all $x\in \tilde{C}_{R/2}$, $\xi\in V_I$. The result follows from here.
\end{proof}

The computation of $q$-Borel transform on polynomials can also be described.

\begin{prop}\label{prop336}
Let $\rho>0$, $d\in\R$ and let $f:D(0,\rho)\cup S_d\to\C$ be given as in Proposition~\ref{prop191}. Let 
$$\varphi(T)=\mathcal{L}_{q;k}^{d}(f(z))(T),\quad T\in\tilde{D}_R.$$
Then, the next identity holds:
$$\mathcal{B}_{q;k}^d(\varphi)(\xi)=\xi,\qquad \xi\in V_{I},$$
for all $\xi \in V_{I}$, for some unbounded sector $V_{I}$ as described in Proposition~\ref{prop302}.
\end{prop}
\begin{proof}
One can assume, without loss of generality, that $d=0$. Observe form the definition of $q$-Borel and $q$-Laplace transform that 
$$\mathcal{B}_{q;k}^{d}(\varphi)(\xi)=-\frac{k i}{2\pi\log(q)}\left[A_1+A_2\right],$$
where
$$A_1:=\int_{\tilde{C}^{+}_{R/2}}q^{\frac{k}{2}\frac{\log(x/\xi)}{\log(q)}\left(\frac{\log(x/\xi)}{\log(q)}-\frac{1}{k}\right)}\int_{L_{d^{+}}}q^{-\frac{k}{2}\frac{\log(x/u)}{\log(q)}\left(\frac{\log(x/u)}{\log(q)}-\frac{1}{k}\right)}f(u)\frac{du}{u}\frac{dx}{x},$$
and
$$A_2:=\int_{\tilde{C}^{-}_{R/2}}q^{\frac{k}{2}\frac{\log(x/\xi)}{\log(q)}\left(\frac{\log(x/\xi)}{\log(q)}-\frac{1}{k}\right)}\int_{L_{d^{-}}}q^{-\frac{k}{2}\frac{\log(x/u)}{\log(q)}\left(\frac{\log(x/u)}{\log(q)}-\frac{1}{k}\right)}f(u)\frac{du}{u}\frac{dx}{x},$$
where $L_{d^{\pm}}=\R_{+}e^{i d_{\pm}}$, for some small $d_{+}>0>d_{-}$, and $\tilde{C}^+_{R/2}=\left\{\frac{R}{2}e^{i\theta}:\theta\ge0\right\}$, $\tilde{C}^-_{R/2}=\left\{\frac{R}{2}e^{i\theta}:\theta\le 0\right\}$.

Classical estimates and the application of Fubini Theorem yields
\begin{multline*}
A_1=\int_{L_{d^{+}}}\int_{\tilde{C}^{+}_{R/2}}e^{\frac{1}{2}\log\left(\frac{\xi}{u}\right)}e^{\frac{k}{2\log(q)}[\log^2(\xi)-\log^2(u)]}f(u)e^{\frac{k}{\log(q)}\log(x)\log(u/\xi)}\frac{dx}{x}\frac{du}{u}\\
=\int_{L_{d^{+}}}\left(\frac{\xi}{u}\right)^{1/2}e^{\frac{k}{2\log(q)}[\log^2(\xi)-\log^2(u)]}f(u)B_1\frac{du}{u},
\end{multline*}
where
$$B_1=\int_{\tilde{C}^{+}_{R/2}}e^{\frac{k}{\log(q)}\log(x)\log(u/\xi)}\frac{dx}{x}
.$$
Taking into account that $x=R/2e^{it}$ and $\hbox{arg}(u/\xi)=\hbox{arg}(u)-\hbox{arg}(\xi)=d^{+}-\hbox{arg}(\xi)$, we derive that $\hbox{arg}(u/\xi)>0$ for all $\xi\in V_{I}$. Recall that $\hbox{arg}(\xi)$ is close to $d=0$ and $d^{+}>0$. Therefore,
$$B_1=ie^{\frac{k}{\log(q)}\log(R/2)\log(u/\xi)}\int_0^{\infty}e^{i\frac{k}{\log(q)}\log(u/\xi)t}dt=\frac{-i}{\frac{k}{\log(q)i\log(u/\xi)}}e^{\frac{k}{\log(q)}\log(R/2)\log(u/\xi)}.$$
As a result,
$$A_1=-\int_{L_{d^+}}\left(\frac{\xi}{u}\right)^{1/2}e^{\frac{k}{2\log(q)}[\log^2(\xi)-\log^2(u)]}f(u)e^{\frac{k}{\log(q)}\log(k/2)\log(u/\xi)}\frac{1}{\frac{k}{\log(q)}\log(u/\xi)}\frac{du}{u}.$$
In a similar manner, one arrives at

$$A_2=\int_{L_{d^{-}}}\left(\frac{\xi}{u}\right)^{1/2}e^{\frac{k}{2\log(q)}[\log^2(\xi)-\log^2(u)]}f(u)B_2\frac{du}{u},$$
with
$$B_2=\int_{\tilde{C}^{-}_{R/2}}e^{\frac{k}{\log(q)}\log(x)\log(u/\xi)}\frac{dx}{x}
.$$
Since $x=R/2e^{it}$ and $\hbox{arg}(u/\xi)=\hbox{arg}(u)-\hbox{arg}(\xi)=d^{-}-\hbox{arg}(\xi)$, we derive that $\hbox{arg}(u/\xi)<0$ for all $\xi\in V_{I}$. Therefore,
$$B_2=ie^{\frac{k}{\log(q)}\log(R/2)\log(u/\xi)}\int_{-\infty}^{0}e^{i\frac{k}{\log(q)}\log(u/\xi)t}dt=\frac{i}{\frac{k}{\log(q)i\log(u/\xi)}}e^{\frac{k}{\log(q)}\log(R/2)\log(u/\xi)}.$$
As a result,
$$A_2=\int_{L_{d^-}}\left(\frac{\xi}{u}\right)^{1/2}e^{\frac{k}{2\log(q)}[\log^2(\xi)-\log^2(u)]}f(u)e^{\frac{k}{\log(q)}\log(k/2)\log(u/\xi)}\frac{1}{\frac{k}{\log(q)}\log(u/\xi)}\frac{du}{u}.$$

We conclude that 
$$\mathcal{B}_{q;k}^{d}(\varphi)(\xi)=\frac{1}{2\pi i}\left[\int_{L_{d^{-}}}A(\xi,u)du+\int_{L_{d^{+}}}A(\xi,u)du\right],$$
with 
$$u\mapsto A(\xi,u)=\left(\frac{\xi}{u}\right)^{1/2}e^{\frac{k}{2\log(q)}[\log^2(\xi)-\log^2(u)]}f(u)e^{\frac{k}{\log(q)}\log(R/2)\log(u/\xi)}\frac{1}{\log(u/\xi)u}.$$
For all $\xi$, the map $A(\xi,u)$ is holomorphic on $U$ if $0<R<1$, and it is clearly bounded for $u\in U$, as $u\to 0$. A geometric configuration is illustrated in Figure~\ref{fig2}.

\begin{figure}
	\centering
\includegraphics[width=0.45\textwidth]{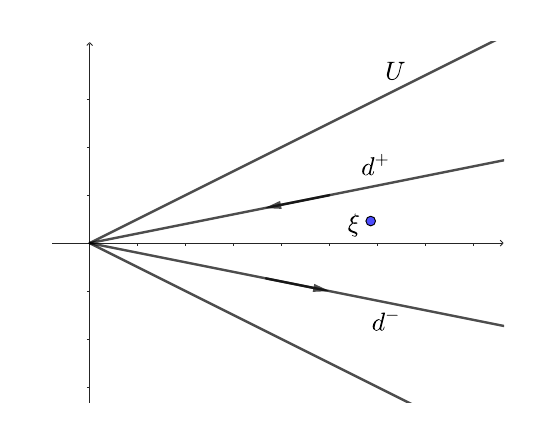}
\caption{Geometric configuration of paths}
\label{fig2}
\end{figure}

We apply the residues formula leading to
$$\displaystyle\mathcal{B}_{q;k}^{d}(\varphi)(\xi)=\hbox{Res}_{u=\xi}(u\mapsto A(\xi,u)),$$
since $u=\xi$ is the only singularity of $u\mapsto A(\xi,u)$ enclosed by the oriented path $L_{d^{-}}-L_{d^{+}}$.

The following classical result can be found in Chapter VI,~\cite{lang}.
\begin{lemma}
Let $U\subseteq\C$ be an open set. The following statements hold.
\begin{enumerate}
\item Let $f:U\setminus\{z_0\}\to\C$ be a holomorphic map with a simple pole at $z_0$ and $g:U\to\C$ be a holomorphic map. Then,
$$\displaystyle\hbox{Res}_{z=z_0}(fg)=g(z_0)\hbox{Res}_{z=z_0}(f).$$
\item Let $f:U\to\C$ be a holomorphic map. Assume that $f(z_0)=0$ and $f'(z_0)\neq 0$. Then,
$$\displaystyle\hbox{Res}_{z=z_0}\left(\frac{1}{f}\right)=\frac{1}{f'(z_0)}.$$
\end{enumerate}
\end{lemma}

In view of the previous Lemma, we arrive at
\begin{multline*}
\displaystyle\hbox{Res}_{u=\xi}(u\mapsto A(\xi,u))=\left(\frac{\xi}{\xi}\right)^{1/2}e^{\frac{k}{2\log(q)}[\log^2(\xi)-\log^2(\xi)]}f(\xi)e^{\frac{k}{\log(q)}\log(R/2)\log(\xi/\xi)}\frac{1}{\xi}\\
\displaystyle\times \hbox{Res}_{u=\xi}\left(\frac{1}{\log(u)-\log(\xi)}\right)=f(\xi).
\end{multline*}

The result follows from here.
\end{proof}

As a result, the next corollary is deduced.

\begin{corol}\label{coro467}
Let $\varphi(x)=x^n$, for $n\ge0$ an integer number. Then,
$$\mathcal{B}^{d}_{q;k}(\varphi)(\xi)=\frac{1}{q^{\frac{n(n-1)}{2k}}}\xi^n.$$
\end{corol}
\begin{proof} 
In view of Proposition~\ref{prop254} and Proposition~\ref{prop336} we have
$$\mathcal{B}^{d}_{q;k}(\mathcal{L}_{q;k}^d(x^n)(T))(\xi)=\xi$$
if and only if 
$$\mathcal{B}^{d}_{q;k}(q^{\frac{n(n-1)}{2k}}T^n)(\xi)=\xi.$$
As a consequence of linearity of the Borel mapping, we arrive at
$$q^{\frac{n(n-1)}{2k}}\mathcal{B}^{d}_{q;k}(T^n)(\xi)=\xi,$$
which yields the result.
\end{proof}

\subsection{Review on $G_q$-summable formal power series. Main properties}\label{sec56}

In this section, we recall the notion of $G_q$-summability of order $k\ge 1$ introduced in~\cite{zh99,mazh}.

\begin{defin}
Let $(\mathbb{E},\left\|\cdot\right\|_{\mathbb{E}})$ be a complex Banach space. Let $k\ge1$ be an integer, and let $q>1$. A formal power series 
$$\hat{U}(T)=\sum_{n\ge 1}a_nT^n\in\mathbb{E}[[T]]$$
is said to be $G_q$-summable of order $k$ with respect to $T$ in the direction $d\in\R$ if the following properties hold:
\begin{itemize}
\item The formal $q$-Borel transform of order $k$ of $\hat{U}(T)$, defined by
$$\hat{\mathcal{B}}_{q;k}(\hat{U})(\xi)=\sum_{n=1}^{\infty}\frac{a_n}{q^{\frac{n(n-1)}{2k}}}\xi^n\in\xi\mathbb{E}[[\xi]]$$
is convergent on some disc $D(0,\rho)$, for some $\rho>0$.
\item The function defining such convergent series can be analytically continued to a function, say $W_k(\xi)$, defined on an unbounded sector of bisecting direction $d$, say
$$S_{d,\delta}=\left\{\xi\in\C^{\star}: |d-\hbox{arg}(\xi)|<\delta\right\}$$
for some $\delta>0$. Moreover, there exist $K,\alpha>0$ and $\Delta>1$ such that 
$$\left\|W_k(\xi)\right\|_{\mathbb{E}}\le K|\xi|\exp\left(\frac{k}{2\log(q)}\log^2(|\xi|+\Delta)+\alpha\log(|\xi|+\Delta)\right),$$
for all $\xi\in S_{d,\delta}$. 
\end{itemize}
\end{defin}

Given a $G_q$-summable formal power series of order $k$ in a direction $d\in\R$, $\hat{U}(T)\in T\mathbb{E}[[T]]$, the so-called $G_q$-sum of $\hat{U}(T)$ along direction $d$ is given by the $q$-Laplace transform
$$U_d(T)=\mathcal{L}_{q;k}^d(W_k)(T)=\pi_{q,k}\int_{L_d}\Theta_k\left(\frac{T}{\tau}\right)W_k(\tau)\frac{d\tau}{\tau}.$$

According to Proposition~\ref{prop191}, $U_d(T)$ turns out to be analytic on the universal covering $\tilde{D}(0,R)$ of the disc $D(0,R)$, for some $R>0$. In addition to this, Proposition~\ref{prop325} guarantees that $U_d(T)$ admits $\hat{U}(T)$ as its asymptotic expansion of $q$-Gevrey type of order $k$ on any sectorial domain $\tilde{S}_d=\{\rho e^{i\theta}:0\le \rho< R,\theta\in I\}$, for any bounded segment $I\subseteq\R$, i.e. there exist $K,A>0$ such that
$$\left\|U_d(T)-\sum_{n=1}^{N-1}a_nT^n\right\|_{\mathbb{E}}\le KA^Nq^{\frac{N^2}{2k}}|T|^{N},$$
for all $N\ge2$ and $T\in\tilde{S}_R$.

\vspace{0.3cm}

\textbf{Aknowledgements:} The first author is partially supported by Laboratoire Paul Painlev\'e de Math\'ematiques, Universit\'e de Lille, by a grant Cross Disciplinary Projects $\hbox{C}^2\hbox{EMPI}$. Both authors are partially supported by the project PID2022-139631NB-I00 of Ministerio de Ciencia e Innovaci\'on, Spain.

\end{document}